\documentclass{amsart}

\usepackage[T1]{fontenc}
\usepackage{lmodern}
\usepackage{microtype}
\usepackage{amssymb,amsmath,amsthm}
\usepackage{setspace}
\usepackage{booktabs}
\usepackage{needspace}
\usepackage{hyperref}
\usepackage[capitalize,noabbrev]{cleveref}

\allowdisplaybreaks
\hypersetup{
  colorlinks=true,
  linkcolor=blue,
  citecolor=blue,
  urlcolor=blue,
  pdftitle={On Grassmann 4-designs},
  pdfauthor={Mohammad Samy Baladram, Ryutaro Misawa, and Akihiro Munemasa},
  pdfsubject={Grassmannian designs and association schemes},
  pdfkeywords={Grassmannian design, association scheme, spherical embedding}
}

\theoremstyle{definition}
\newtheorem{theorem}{Theorem}[section]
\newtheorem{proposition}[theorem]{Proposition}
\newtheorem{lemma}[theorem]{Lemma}
\newtheorem{definition}[theorem]{Definition}
\theoremstyle{remark}
\newtheorem{remark}[theorem]{Remark}
\numberwithin{equation}{section}

\DeclareMathOperator{\rank}{rank}
\DeclareMathOperator{\tr}{tr}
\DeclareMathOperator{\Sym}{Sym}

\title[On Grassmann $4$-designs]{On Grassmann $4$-designs}

\author[M. S. Baladram]{Mohammad Samy Baladram}
\address{Undergraduate Admissions Division, Admissions Center, Tohoku University, Sendai, Japan}
\email{samy.baladram@tohoku.ac.jp}

\author[R. Misawa]{Ryutaro Misawa}
\address{Graduate School of Information Sciences, Tohoku University, Sendai 980-8579, Japan}
\email{misawa.ryutaro.q2@dc.tohoku.ac.jp}

\author[A. Munemasa]{Akihiro Munemasa}
\address{Tohoku University and China University of Geosciences}
\email{munemasa@tohoku.ac.jp}


\begin{document}
\raggedbottom

\begin{abstract}
We prove that there is a unique tight $4$-design in the real Grassmannian $G_{2,4}$, up to the action of $O(4)$.
Using symmetric association schemes, we also construct an $N$-point $4$-design in $G_{2,4}$ for every even integer $N\geq10$.
For symmetric association schemes of order less than $21$, we classify the pairs of rank-three primitive idempotents defining injective maps into $G_{2,4}$; their images are Grassmannian $4$-designs precisely for four schemes of order $16$.
\end{abstract}

\maketitle

\section{Introduction}

A Grassmannian design is a configuration that reproduces the averages of polynomial functions of degree up to a prescribed bound on a Grassmannian as equi-weighted averages over finitely many subspaces.
For integers $1\leq k<d$, let $G_{k,d}:=\{V\leq\mathbb R^d:\dim V=k\}$ denote the real Grassmannian, and take averages on this space with respect to the orthogonally invariant probability measure.
A spherical design is a configuration that reproduces the spherical averages of polynomials of degree up to a prescribed bound as equi-weighted averages over a finite point set \cite[Section~5]{DelsarteGoethalsSeidel1977}.
Spherical designs have been extensively studied with regard to existence, explicit constructions, lower bounds on cardinality, and the structure of configurations of minimum cardinality \cite{DelsarteGoethalsSeidel1977,deLaHarpePache2005}.
Against this background in spherical design theory, Bachoc, Coulangeon, and Nebe \cite{BachocCoulangeonNebe2002} introduced Grassmannian designs as an analogue for configurations of subspaces.

Like spherical designs, Grassmannian designs satisfy lower bounds on their cardinality in terms of dimension and strength, and designs attaining these bounds are called tight \cite[Theorem~8]{BachocBannaiCoulangeon2004}.
In the spherical case, the tight designs known to exist are unique up to orthogonal transformations for each dimension and strength \cite[Section~1.17]{deLaHarpePache2005}.
For Grassmannian designs, however, uniqueness need not hold.
Indeed, the construction of Calderbank, Hardin, Rains, Shor, and Sloane \cite[Section~5]{CalderbankHardinRainsShorSloane1999} yields noncongruent tight $4$-designs with $28$ points in $G_{3,7}$ \cite[Section~4.2]{BachocBannaiCoulangeon2004}.
Thus, in addition to existence, determining when uniqueness holds is a problem for tight Grassmannian designs.

Grassmannian designs have been studied in connection with lattice theory and the representation theory of finite groups.
Bachoc, Coulangeon, and Nebe \cite{BachocCoulangeonNebe2002} examined conditions under which configurations of subspaces arising from lattices are Grassmannian designs and related this property to the extremality of lattices.
Harmonic analysis on Grassmannians has also been used to obtain criteria for designs and lower bounds on their cardinality \cite{BachocBannaiCoulangeon2004,Bachoc2005}.
In this way, Grassmannian designs have been studied both through the construction of configurations of subspaces and through their geometric and algebraic structures.

In particular, a Grassmannian $2$-design is equivalent to an equi-weighted tight fusion frame consisting of subspaces of the same dimension \cite[Section~5]{BachocEhler2013}, \cite[Remark~2.18]{Misawa2026}.
This is a configuration for which the sum of the orthogonal projections onto the subspaces is a scalar multiple of the identity operator, and conditions for the existence of such configurations and their explicit constructions have also been studied from the viewpoint of frame theory \cite[Sections~3--4]{CalderbankCasazzaHeineckeKutyniokPezeshki2011}, \cite[Section~6]{BachocEhler2013}.
Here, tightness of a fusion frame refers to a condition on the sum of the projection operators and is distinct from tightness of a design, which means attaining the standard cardinality bound.

On the other hand, if $2\leq k\leq d-2$, there are no tight Grassmannian $2$-designs in $G_{k,d}$.
Indeed, applying the standard bound above at strength $2$ gives a lower bound of $1$, but a Grassmannian $2$-design is equivalent to an equi-weighted tight fusion frame, and a single subspace cannot satisfy this condition \cite[Proposition~4.2]{BachocCoulangeonNebe2002}, \cite[Theorem~8]{BachocBannaiCoulangeon2004}, \cite[Remark~2.18]{Misawa2026}.
Thus, $4$ is the next strength to consider when studying the existence and uniqueness of tight Grassmannian designs.

In this paper, we study the construction and uniqueness of Grassmannian $4$-designs in $G_{2,4}$, a basic nontrivial case.
With $S^2:=\{u\in\mathbb R^3:\|u\|=1\}$, this space can be represented as
\[
G_{2,4}\cong(S^2\times S^2)/{\sim}.
\]
Here, $\sim$ is the equivalence relation identifying $(\ell,r)$ with $(-\ell,-r)$ \cite[Section~4]{ConwayHardinSloane1996}.
This representation allows configurations of planes to be studied through configurations on two spheres.
Finite configurations in $G_{2,4}$ have also been studied from the perspectives of packing and energy minimization \cite{ConwayHardinSloane1996,EtayoLopezGomez2025}.
In this paper, we use this relationship with spheres to study the structure and construction of Grassmannian designs.

We first consider the uniqueness of tight Grassmannian $4$-designs.
The definition of a Grassmannian $4$-design in $G_{2,4}$ is given in \cref{def:grassmann-design}.
Such a design has at least $10$ elements, and a configuration of $10$ planes attaining this bound was constructed by Conway, Hardin, and Sloane \cite[Sections~4--5]{ConwayHardinSloane1996}, \cite{BachocBannaiCoulangeon2004}, \cite[Example~3.7(1) and Remark~5.6]{BachocEhler2013}.
We prove that every tight Grassmannian $4$-design can be mapped onto this configuration by an orthogonal transformation of $\mathbb R^4$ (\cref{thm:tight-uniqueness}).

We next consider the construction of Grassmannian $4$-designs.
One existing method for constructing Grassmannian designs uses orbits of finite subgroups of the orthogonal group.
Bachoc, Coulangeon, and Nebe \cite[Theorem~4.5]{BachocCoulangeonNebe2002} and Bachoc \cite[Section~4]{Bachoc2005} give representation-theoretic conditions under which finite group orbits are Grassmannian designs.
This method extends the construction of spherical designs by group orbits, together with the corresponding criteria, to Grassmannians.
Spherical embeddings of symmetric association schemes also play an important role in the construction of spherical designs \cite[Section~9]{DelsarteGoethalsSeidel1977}, \cite{Suda2011,LansdownMartinMunemasaSudaTanaka2026}.
In this paper, we consider a Grassmannian analogue of this construction method.

Specifically, from a symmetric association scheme, we obtain two configurations on $S^2$ indexed by the same set.
We combine them using the above representation of planes to construct a map into $G_{2,4}$.
When this map is injective, we give necessary and sufficient conditions for its image to be a Grassmannian $4$-design (\cref{thm:association-scheme-construction}).
Applying this to the association scheme of an even cycle, we construct a Grassmannian $4$-design with $N$ elements for every even $N\geq10$ (\cref{thm:even-cycle}).
For $N=10$, this recovers the known configuration of Conway, Hardin, and Sloane.

We further restrict attention to pairs of rank-three primitive idempotents in the construction using idempotents formulated in \cref{sec:association-scheme-construction}.
For symmetric association schemes on at most $20$ points, we classify all such pairs defining injective maps into $G_{2,4}$ (\cref{thm:small-order-classification}).
Pairs whose associated maps have Grassmannian $4$-designs as their images arise only from four schemes on $16$ points.
The idempotents used in the construction from even cycles are not primitive, so this construction is not included in the classification.
This classification uses the classification of spherical embeddings by Bannai and Zhao \cite{BannaiZhao2020} and the complete lists of association schemes of the relevant orders \cite{Nomiyama1995,Hirasaka1997,HanakiMiyamoto1998a,HanakiMiyamoto1998b,HanakiMiyamoto2003}.

Section~\ref{sec:gram-criterion} gives a criterion in terms of Gram matrices, and \cref{sec:association-scheme-construction} applies it to association schemes.
Section~\ref{sec:even-cycles} presents the construction from even cycles, and \cref{sec:tight-uniqueness} proves the uniqueness of tight designs.
Section~\ref{sec:primitive-idempotents} treats primitive idempotents and gives the classification for symmetric association schemes on at most $20$ points.

\section{A Gram-matrix criterion}
\label{sec:gram-criterion}

For real vectors, write $\langle x,y\rangle:=x^{\mathsf T}y$ and $\lVert x\rVert:=\sqrt{\langle x,x\rangle}$.
For $n\geq1$, let $I_n$ and $J_n$ be the identity and all-ones matrices of order $n$, respectively.
Write $\Sym(n)$ for the space of real symmetric $n\times n$ matrices, $O(n)$ for the orthogonal group on $\mathbb R^n$, and $SO(n)$ for its subgroup of determinant $1$.
For $V\in G_{2,4}$, let $\Pi_V:\mathbb R^4\to\mathbb R^4$ be the orthogonal projection with image $V$, identified with its matrix in the standard basis.
Let $\sigma$ be the unique $O(4)$-invariant probability measure on $G_{2,4}$.

\Needspace{8\baselineskip}
\begin{definition}
\label[definition]{def:grassmann-design}
A finite nonempty set $D\subset G_{2,4}$ is a \emph{Grassmannian $4$-design} if
\[
  \frac1{|D|}\sum_{V\in D}p(\Pi_V)
  =\int_{G_{2,4}}p(\Pi_V)\,d\sigma(V)
\]
for every real polynomial function $p:\Sym(4)\to\mathbb R$ of total degree at most $2$ in the matrix entries.
It is \emph{tight} if $|D|=10$.
\end{definition}

\Cref{def:grassmann-design} agrees with \cite[Definition~4.1 and Proposition~4.2]{BachocCoulangeonNebe2002}.
To see this, consider real polynomials $f$ of total degree at most $4$ in the eight matrix entries of $U\in\mathbb R^{4\times2}$, satisfying
\[
  f(UQ)=f(U)\qquad(U\in\mathbb R^{4\times2},\ Q\in O(2)).
\]
By \cite[Theorem~2.9.A]{Weyl1946}, applied to the four rows of $U$, these are precisely the polynomials of the form $f(U)=p(UU^{\mathsf T})$ with $\deg p\leq2$.
For each $V\in G_{2,4}$, we can choose $U\in\mathbb R^{4\times2}$ whose columns form an orthonormal basis of $V$, so that $\Pi_V=UU^{\mathsf T}$.
Then $f(U)=p(\Pi_V)$, independently of this choice of $U$, so the two definitions require the same averaging identities.

To express the design condition in terms of Gram matrices, identify $\mathbb R^4$ with the quaternion algebra $\mathbb H$, and $\mathbb R^3$ with its purely imaginary subspace.
For $\ell,r\in S^2$, put
\[
  [\ell,r]:=\{x\in\mathbb H:-\ell xr=x\}.
\]
By \cite[Theorem~4.1]{ConwayHardinSloane1996}, $[\ell,r]\in G_{2,4}$ and every element of $G_{2,4}$ has this form.
For $V=[\ell,r]$ and $W=[\ell',r']$ with $\ell,r,\ell',r'\in S^2$,
\[
  V=W
  \quad\Longleftrightarrow\quad
  (\ell',r')\in\{(\ell,r),(-\ell,-r)\},
\]
and
\begin{equation}
  \tr(\Pi_V\Pi_W)=1+\langle\ell,\ell'\rangle\langle r,r'\rangle.
  \label{eq:projection-inner-product}
\end{equation}
Equation \eqref{eq:projection-inner-product} follows from \cite[Equation~(4.1) and Section~5]{ConwayHardinSloane1996}.

To apply the design criterion in \cite[Theorem~6.1]{Bachoc2005}, we compute $\int_{G_{2,4}}(\tr(\Pi_V\Pi_W))^2\,d\sigma(W)$ for fixed $V\in G_{2,4}$.
Let $\mu$ be the uniform probability measure on $S^2$.
For unit quaternions $q,s$, the map $x\mapsto qxs^{-1}$ is an element of $SO(4)$, and every element of $SO(4)$ has this form \cite[Section~4]{ConwayHardinSloane1996}.
If $x\in[\ell,r]$, then
\[
  -(q\ell q^{-1})(qxs^{-1})(srs^{-1})
  =q(-\ell xr)s^{-1}
  =qxs^{-1}.
\]
The inverse map $x\mapsto q^{-1}xs$ sends $[q\ell q^{-1},srs^{-1}]$ into $[\ell,r]$, so $x\mapsto qxs^{-1}$ sends $[\ell,r]$ onto $[q\ell q^{-1},srs^{-1}]$.
Conjugation by a unit quaternion fixes the real axis and belongs to $SO(4)$, so its restriction to the purely imaginary subspace is a rotation of $\mathbb R^3$.
Thus conjugations by $q$ and $s$ preserve $\mu\times\mu$.
Conversely, extending any rotation of $\mathbb R^3$ by the identity on the real axis gives an element of $SO(4)$ of the form $x\mapsto qxs^{-1}$; evaluating it at $1$ gives $qs^{-1}=1$, hence $q=s$.
Quaternionic conjugation $x\mapsto\bar x$ sends $[\ell,r]$ onto $[r,\ell]$, since $\bar\ell=-\ell$ and $\bar r=-r$ give
\[
  -r\bar x\ell
  =\overline{-\ell xr}
  =\bar x
  \qquad(x\in[\ell,r]),
\]
and applying conjugation twice gives the inverse inclusion.
The interchange $(\ell,r)\mapsto(r,\ell)$ preserves $\mu\times\mu$, and quaternionic conjugation has determinant $-1$, because it is the identity on the real axis and minus the identity on the purely imaginary subspace.
Together with $SO(4)$, it generates $O(4)$, so the measure induced by $\mu\times\mu$ under $(\ell,r)\mapsto[\ell,r]$ is $O(4)$-invariant and equals $\sigma$.
For fixed $a\in S^2$, the average of $\langle a,u\rangle$ with respect to $\mu$ is $0$, and that of $\langle a,u\rangle^2$ is $1/3$ by rotational invariance and $\lVert u\rVert^2=1$.
Thus, for $V=[\ell,r]\in G_{2,4}$,
\begin{equation}
  \begin{aligned}
  \int_{G_{2,4}}\bigl(\tr(\Pi_V\Pi_W)\bigr)^2\,d\sigma(W)
  &=\int_{S^2}\int_{S^2}
    \bigl(1+\langle\ell,u\rangle\langle r,v\rangle\bigr)^2\,d\mu(u)\,d\mu(v)=\frac{10}{9}.
  \end{aligned}
  \label{eq:invariant-second-moment}
\end{equation}

For real matrices, $\circ$ denotes the Hadamard product of matrices of the same size, $\lVert\cdot\rVert_F$ the Frobenius norm, and $\otimes$ the Kronecker product.

Let $X$ be a finite nonempty set, and let $D=\{[\ell_x,r_x]:x\in X\}\subset G_{2,4}$ consist of $N:=|X|$ distinct elements, with $\ell_x,r_x\in S^2$.
Let $L,R$ be the Gram matrices of $(\ell_x)_{x\in X}$ and $(r_x)_{x\in X}$, respectively.
For $x\in X$, put $z_x:=\ell_x\otimes r_x\in\mathbb R^9$.
Then $L\circ R$ is the Gram matrix of the unit vectors $z_x$.
Replacing $(\ell_x,r_x)$ by $(-\ell_x,-r_x)$ changes neither $z_x$ nor $L\circ R$.

\begin{lemma}
\label[lemma]{lem:G24-moment-criterion}
The set $D$ is a Grassmannian $4$-design if and only if
\begin{equation}
  \sum_{x,y\in X}(L\circ R)_{xy}=0
  \label{eq:first-moment}
\end{equation}
and
\begin{equation}
  \lVert L\circ R\rVert_F^2=\frac{N^2}{9}.
  \label{eq:second-product-moment}
\end{equation}
Equivalently,
\begin{equation}
  \sum_{x\in X}z_x=0,
  \qquad
  \sum_{x\in X}z_xz_x^{\mathsf T}=\frac N9I_9.
  \label{eq:tensor-moments}
\end{equation}
In particular, $N\geq10$ whenever $D$ is a Grassmannian $4$-design.
\end{lemma}

\begin{proof}
By \cite[Theorem~6.1]{Bachoc2005} and \eqref{eq:invariant-second-moment}, $D$ is a Grassmannian $4$-design if and only if the average of $(\tr(\Pi_V\Pi_W))^2$ over $D\times D$ is $10/9$.
Put $T:=\sum_{x\in X}z_xz_x^{\mathsf T}$, so that $\tr(T)=N$.
Equation \eqref{eq:projection-inner-product} gives
\begin{equation}
  \frac1{N^2}\sum_{V,W\in D}\bigl(\tr(\Pi_V\Pi_W)\bigr)^2
  =1+\frac2{N^2}\left\lVert\sum_{x\in X}z_x\right\rVert^2+\frac1{N^2}\tr(T^2).
  \label{eq:trace-square-sum}
\end{equation}
Moreover, $\tr(T)=N$ implies
\begin{equation}
  \tr(T^2)-\frac{N^2}{9}
  =\left\lVert T-\frac N9I_9\right\rVert_F^2\geq0.
  \label{eq:tight-frame-bound}
\end{equation}
Thus the average in \eqref{eq:trace-square-sum} is $10/9$ if and only if \eqref{eq:tensor-moments} holds.
The identities
\[
  \sum_{x,y\in X}(L\circ R)_{xy}
  =\left\lVert\sum_{x\in X}z_x\right\rVert^2,
  \qquad
  \lVert L\circ R\rVert_F^2=\tr(T^2),
\]
together with \eqref{eq:tight-frame-bound}, show that \eqref{eq:tensor-moments} is equivalent to \eqref{eq:first-moment} and \eqref{eq:second-product-moment}.
Finally, $T=(N/9)I_9$ shows that the vectors $z_x$ span $\mathbb R^9$, while $\sum_{x\in X}z_x=0$ expresses any one of them as a linear combination of the others.
The remaining $N-1$ vectors therefore still span $\mathbb R^9$, so $N-1\geq9$.
\end{proof}

\section{Construction from association schemes}
\label{sec:association-scheme-construction}

We construct Gram matrices from the commuting projections in the Bose--Mesner algebra of a symmetric association scheme and apply \cref{lem:G24-moment-criterion}.
We use the standard definition and algebraic properties of association schemes from \cite[Section~2]{MartinTanaka2009}.

\begin{definition}
\label[definition]{def:symmetric-association-scheme}
Let $X$ be a finite nonempty set, put $N:=|X|$, and let $\{R_i\}_{i=0}^s$ be a partition of $X^2$ into nonempty relations.
Define the adjacency matrix $A_i$ by $(A_i)_{xy}:=1$ if $(x,y)\in R_i$ and $(A_i)_{xy}:=0$ otherwise.
Then $(X,\{R_i\}_{i=0}^s)$ is a \emph{symmetric association scheme} if
\[
  A_0=I_N,\qquad A_i^{\mathsf T}=A_i,\qquad
  A_iA_j\in\operatorname{span}_{\mathbb R}\{A_0,\ldots,A_s\}
  \qquad(0\leq i,j\leq s).
\]
Its \emph{Bose--Mesner algebra} is
\[
  \mathcal A:=\operatorname{span}_{\mathbb R}\{A_0,\ldots,A_s\}.
\]
\end{definition}

For $A,B\in\mathcal A$, closure under multiplication and symmetry give $AB=(AB)^{\mathsf T}=BA$.
Thus simultaneous diagonalization gives primitive idempotents $E_0,\ldots,E_s$, with $E_0=J_N/N$, which form a basis of $\mathcal A$ and satisfy
\[
  E_jE_k=\delta_{jk}E_j\quad(0\leq j,k\leq s),\qquad
  \sum_{j=0}^sE_j=I_N,
\]
where $\delta_{jk}$ is the Kronecker delta.
Put $m_j:=\rank(E_j)$, the multiplicity of $E_j$.
Since $A_0$ is the only adjacency matrix with nonzero diagonal entries, every matrix in $\mathcal A$ has constant diagonal entries.
An idempotent $F\in\mathcal A$ is an orthogonal projection, so its diagonal entries are $\rank(F)/N$.
For an idempotent $F\in\mathcal A$ of rank $3$, the matrix $NF/3$ is therefore positive semidefinite of rank $3$ with diagonal entries $1$.
Thus there is a map $\varphi_F:X\to S^2$ that satisfies
\begin{equation}
  \langle\varphi_F(x),\varphi_F(y)\rangle=\frac N3F_{xy}
  \qquad(x,y\in X).
  \label{eq:spherical-embedding-gram}
\end{equation}
When $F$ is primitive, $\varphi_F$ is its \emph{spherical embedding} \cite{BannaiZhao2020}.
We call $\varphi_F$ \emph{faithful} if $\varphi_F$ is injective.

For idempotents $F,G\in\mathcal A$ of rank $3$, choose maps $\varphi_F,\varphi_G:X\to S^2$ with Gram matrices $NF/3$ and $NG/3$, respectively, and define
\begin{equation}
  \Phi_{F,G}:X\longrightarrow G_{2,4},
  \qquad x\longmapsto[\varphi_F(x),\varphi_G(x)].
  \label{eq:scheme-map}
\end{equation}

\begin{lemma}
\label[lemma]{lem:faithfulness}
The map $\Phi_{F,G}$ is injective if and only if
\[
  \left(\frac N3F_{xy},\frac N3G_{xy}\right)\notin\{(1,1),(-1,-1)\}
  \qquad(x,y\in X,\ x\neq y).
\]
\end{lemma}

\begin{proof}
By the representation of planes in \cref{sec:gram-criterion} and \eqref{eq:spherical-embedding-gram},
\[
  \begin{aligned}
    \Phi_{F,G}(x)=\Phi_{F,G}(y)
    &\iff (\varphi_F(x),\varphi_G(x))
           =\pm(\varphi_F(y),\varphi_G(y))\\
    &\iff \left(\frac N3F_{xy},\frac N3G_{xy}\right)
           \in\{(1,1),(-1,-1)\}.
  \end{aligned}
\]
\end{proof}

\begin{theorem}
\label[theorem]{thm:association-scheme-construction}
Let $F,G\in\mathcal A$ be idempotents of rank $3$, and suppose that $\Phi_{F,G}$ is injective.
Then $\Phi_{F,G}(X)$ is a Grassmannian $4$-design if and only if
\[
  FG=0
  \qquad\text{and}\qquad
  \bigl(N(F\circ G)\bigr)^2=N(F\circ G).
\]
\end{theorem}

\begin{proof}
Let $L:=NF/3$ and $R:=NG/3$.
Since $F$ and $G$ are symmetric,
\begin{equation}
  \sum_{x,y\in X}(L\circ R)_{xy}
  =\frac{N^2}{9}\tr(FG).
  \label{eq:scheme-first-moment}
\end{equation}
Since $FG$ is an idempotent, $\tr(FG)=0$ if and only if $FG=0$.

Suppose first that $\Phi_{F,G}(X)$ is a Grassmannian $4$-design.
Equation \eqref{eq:scheme-first-moment} and \cref{lem:G24-moment-criterion} give $FG=0$.
Let $Z\in\mathbb R^{9\times N}$ have columns $\varphi_F(x)\otimes\varphi_G(x)$, indexed by $x\in X$.
Then $Z^{\mathsf T}Z=L\circ R$, and \eqref{eq:tensor-moments} gives $ZZ^{\mathsf T}=(N/9)I_9$.
Hence
\[
  (L\circ R)^2
  =Z^{\mathsf T}(ZZ^{\mathsf T})Z
  =\frac N9(L\circ R).
\]
Since $N(F\circ G)=(9/N)(L\circ R)$, we obtain $\bigl(N(F\circ G)\bigr)^2=N(F\circ G)$.

Conversely, suppose that $FG=0$ and $N(F\circ G)$ is idempotent.
The diagonal entries of $N(F\circ G)$ are $9/N$, so $\tr(N(F\circ G))=9$.
Since $N(F\circ G)$ is symmetric,
\[
  \lVert L\circ R\rVert_F^2
  =\frac{N^2}{81}\tr\bigl((N(F\circ G))^2\bigr)
  =\frac{N^2}{9}.
\]
Together with \eqref{eq:scheme-first-moment}, this gives \eqref{eq:first-moment} and \eqref{eq:second-product-moment}.
\end{proof}

\section{Construction from even cycles}
\label{sec:even-cycles}

We apply \cref{thm:association-scheme-construction} to the association schemes of cycle graphs.
Let $N\geq4$ be even, put $X:=\mathbb Z/N\mathbb Z$, and define
\[
  R_i:=\{(x,y)\in X^2:x-y=\pm i\}
  \qquad(0\leq i\leq N/2).
\]
These relations form the association scheme of the cycle graph on $X$ \cite{BannaiIto1984}.
For $x,y\in X$, the primitive idempotents $E_0,\ldots,E_{N/2}$ have entries
\begin{equation}
  (E_j)_{xy}=
  \begin{cases}
    1/N,&\text{if }j=0,\\
    (2/N)\cos\bigl(2\pi j(x-y)/N\bigr),&\text{if }1\leq j<N/2,\\
    (-1)^{x-y}/N,&\text{if }j=N/2.
  \end{cases}
  \label{eq:cycle-idempotent-entries}
\end{equation}
These expressions are independent of the integer representatives of $x$ and $y$.
The images of $E_0$ and $E_{N/2}$ are spanned by the constant functions and $x\mapsto(-1)^x$, respectively.
For $1\leq j<N/2$, the image of $E_j$ is spanned by $x\mapsto\cos(2\pi jx/N)$ and $x\mapsto\sin(2\pi jx/N)$.
The orthogonal projections onto these subspaces have entries \eqref{eq:cycle-idempotent-entries}, with $m_0=m_{N/2}=1$ and $m_j=2$ for $1\leq j<N/2$.
Since $E_0=J_N/N$, we have $N(E_0\circ E_j)=E_j$ for $0\leq j\leq N/2$.
For $1\leq i<j$ with $i+j<N/2$, we have
\[
  N(E_i\circ E_j)=E_{i+j}+E_{j-i}.
\]
Multiplying $(E_j)_{xy}$ by $(-1)^{x-y}$ also gives $N(E_j\circ E_{N/2})=E_{N/2-j}$ for $1\leq j<N/2$.

Define unit vectors
\begin{equation}
  \ell_x:=\left(\sqrt{\frac23}\cos\frac{2\pi x}{N},\sqrt{\frac23}\sin\frac{2\pi x}{N},\frac1{\sqrt3}\right)^{\mathsf T},
  \qquad x\in X,
  \label{eq:cycle-left-vectors}
\end{equation}
and
\begin{equation}
  r_x:=\left(\sqrt{\frac23}\cos\frac{4\pi x}{N},-\sqrt{\frac23}\sin\frac{4\pi x}{N},\frac{(-1)^x}{\sqrt3}\right)^{\mathsf T},
  \qquad x\in X.
  \label{eq:cycle-right-vectors}
\end{equation}
\begin{theorem}
\label[theorem]{thm:even-cycle}
For every even integer $N\geq10$, $G_{2,4}$ admits an $N$-point Grassmannian $4$-design.
\end{theorem}

\begin{proof}
Put $F:=E_0+E_1$ and $G:=E_2+E_{N/2}$.
The indices $0,1,2,N/2$ are distinct, so $\rank(F)=\rank(G)=3$ and $FG=0$.
The Gram matrices of the $\ell_x$ and the $r_x$ are $NF/3$ and $NG/3$, respectively.
For distinct $x,y\in X$, we have $-1/3\leq\langle\ell_x,\ell_y\rangle<1$, so \cref{lem:faithfulness} shows that $\Phi_{F,G}$ is injective.

Expanding $N((E_0+E_1)\circ(E_2+E_{N/2}))$ gives
\begin{equation}
  N(F\circ G)=E_1+E_2+E_3+E_{N/2-1}+E_{N/2}.
  \label{eq:cycle-idempotent-product}
\end{equation}
Since $N\geq10$, the indices $1,2,3,N/2-1,N/2$ are distinct, so $N(F\circ G)$ is a sum of mutually orthogonal idempotents and is therefore idempotent.
By \cref{thm:association-scheme-construction}, $\Phi_{F,G}(X)$ is an $N$-point Grassmannian $4$-design.
\end{proof}

For $N=10$, write $u_j:=\ell_j$ and $v_j:=r_j$ for $j\in\mathbb Z/10\mathbb Z$, with $\ell_j,r_j$ as in \eqref{eq:cycle-left-vectors} and \eqref{eq:cycle-right-vectors}, and put
\[
  D_{10}:=\{[u_j,v_j]:j\in\mathbb Z/10\mathbb Z\}.
\]
Let $L$ and $R$ be the Gram matrices of the $u_j$ and the $v_j$, respectively.
Using $L=(10/3)(E_0+E_1)$, $R=(10/3)(E_2+E_5)$, and $E_0=J_{10}/10$, we obtain
\begin{equation}
  L\circ R=\frac{10}{9}\sum_{h=1}^5E_h=\frac{10}{9}I_{10}-\frac19J_{10}.
  \label{eq:D10-Hadamard}
\end{equation}
In particular, $\langle u_j,u_k\rangle\langle v_j,v_k\rangle=-1/9$ for distinct $j,k\in\mathbb Z/10\mathbb Z$.
Moreover, $v_j=(-1)^ju_{3j}$ for $j\in\mathbb Z/10\mathbb Z$.
Thus $D_{10}$ is the $10$-plane design of Conway, Hardin, and Sloane \cite[Section~4]{ConwayHardinSloane1996}.

\begin{remark}
\label[remark]{rem:D10-association-scheme}
For $N=10$, put $V_x:=[\ell_x,r_x]$ for $x\in\mathbb Z/10\mathbb Z$.
For $x,y\in\mathbb Z/10\mathbb Z$, put $\alpha:=\arccos\langle\ell_x,\ell_y\rangle$ and $\beta:=\arccos\langle r_x,r_y\rangle$.
By \cite[Equation~(4.1)]{ConwayHardinSloane1996}, the two eigenvalues of $\Pi_{V_x}\Pi_{V_y}$ that may be nonzero are $\cos^2((\alpha+\beta)/2)$ and $\cos^2((\alpha-\beta)/2)$.
Hence
\begin{align*}
  \tr((\Pi_{V_x}\Pi_{V_y})^2)
  &=\cos^4\frac{\alpha+\beta}{2}+\cos^4\frac{\alpha-\beta}{2}\\
  &=(1+\cos\alpha\cos\beta)^2-\frac12(\cos\alpha+\cos\beta)^2\\
  &=\bigl(1+\langle\ell_x,\ell_y\rangle\langle r_x,r_y\rangle\bigr)^2\\
  &\qquad-\frac12\bigl(\langle\ell_x,\ell_y\rangle+\langle r_x,r_y\rangle\bigr)^2.
\end{align*}
Substituting the inner products obtained from \eqref{eq:cycle-left-vectors} and \eqref{eq:cycle-right-vectors} gives
\[
  \tr((\Pi_{V_x}\Pi_{V_y})^2)=
  \begin{cases}
    2,&\text{if }x-y=0,\\
    64/81,&\text{if }x-y=5,\\
    83/162,&\text{if }x-y\in\{\pm1,\pm3\},\\
    119/162,&\text{if }x-y\in\{\pm2,\pm4\}.
  \end{cases}
\]
Writing $A_i$ for the adjacency matrix of $R_i$, the adjacency matrices for the trace values $2$, $64/81$, $83/162$, and $119/162$, in this order, are
\[
  A_0,\qquad A_5,\qquad A_1+A_3,\qquad A_2+A_4.
\]
Under the bijection $x\mapsto(x\bmod2,x\bmod5)$, these matrices correspond, respectively, to
\[
  \begin{aligned}
    &I_2\otimes I_5,\qquad (J_2-I_2)\otimes I_5,\\
    &(J_2-I_2)\otimes(J_5-I_5),\qquad I_2\otimes(J_5-I_5).
  \end{aligned}
\]
They are the adjacency matrices of the direct product of the one-class association schemes on $2$ and $5$ points \cite{BannaiIto1984}.
This scheme is listed as \texttt{as10[6]} in the catalogue of Hanaki and Miyamoto \cite{HanakiMiyamotoCatalogue}.
\end{remark}

\section{Uniqueness of a tight \texorpdfstring{$4$}{4}-design}
\label{sec:tight-uniqueness}

\begin{theorem}
\label[theorem]{thm:tight-uniqueness}
The design $D_{10}$ is the unique tight $4$-design in $G_{2,4}$, up to the action of $O(4)$.
\end{theorem}

For a tight design $D=\{[\ell_i,r_i]:1\leq i\leq10\}$, we first find, after interchanging $\ell_i$ and $r_i$ for every $i$ if necessary, an invertible $H\in\Sym(3)$ such that $\tr H=0$ and $\ell_i^{\mathsf T}H\ell_i=0$ for every $i$.
The vectors $\ell_i$ lie on the quadric $\{x\in\mathbb R^3:x^{\mathsf T}Hx=0\}$.
The parametrization in \cref{lem:finite-conic} and the characteristic polynomial in \cref{lem:multiplication-matrix} determine the Gram matrix of the $\ell_i$; \eqref{eq:pair-condition} then gives that of the $r_i$.
Finally, we show that an element of $O(4)$ maps $D_{10}$ to $D$.

For a tight design $D=\{[\ell_i,r_i]:1\leq i\leq10\}$, condition \eqref{eq:pair-condition} holds.
Indeed, if they form such a design, \cref{lem:G24-moment-criterion} gives $Z(1,\ldots,1)^{\mathsf T}=0$ and $ZZ^{\mathsf T}=(10/9)I_9$ for the matrix $Z$ with columns $\ell_i\otimes r_i$.
Thus $Z$ has rank $9$ and row space $(1,\ldots,1)^\perp$, so $Z^{\mathsf T}Z=(10/9)I_{10}-(1/9)J_{10}$, which is \eqref{eq:pair-condition}.

\begin{lemma}
\label[lemma]{lem:pair-moments}
Let $\ell_i,r_i\in S^2$ $(1\leq i\leq10)$ satisfy
\begin{equation}
  \langle\ell_i,\ell_j\rangle\langle r_i,r_j\rangle
  =\frac{10\delta_{ij}-1}{9}
  \qquad(1\leq i,j\leq10).
  \label{eq:pair-condition}
\end{equation}
Then
\begin{equation}
  \sum_{i=1}^{10}\ell_i r_i^{\mathsf T}=0,
  \label{eq:pair-center}
\end{equation}
\begin{equation}
  \sum_{i=1}^{10}\ell_i\ell_i^{\mathsf T}
  =\sum_{i=1}^{10}r_i r_i^{\mathsf T}
  =\frac{10}{3}I_3,
  \label{eq:pair-frame}
\end{equation}
and, for all $A,B\in\Sym(3)$,
\begin{equation}
  \sum_{i=1}^{10}(\ell_i^{\mathsf T}A\ell_i)(r_i^{\mathsf T}Br_i)
  =\frac{10}{9}\tr(A)\tr(B).
  \label{eq:pair-quadratic}
\end{equation}
Moreover, $\ell_i\neq\pm\ell_j$ and $r_i\neq\pm r_j$ whenever $i\neq j$.
\end{lemma}

\begin{proof}
Let $Z$ be the $9\times10$ matrix with columns $\ell_i\otimes r_i$.
Equation \eqref{eq:pair-condition} gives
\[
  Z^{\mathsf T}Z=\frac{10}{9}I_{10}-\frac19J_{10}.
\]
The kernel of $Z^{\mathsf T}Z$ is $\mathbb R(1,\ldots,1)^{\mathsf T}$, and its nine nonzero eigenvalues are all $10/9$.
Consequently, $\sum_i\ell_i\otimes r_i=0$, which is \eqref{eq:pair-center}, and $ZZ^{\mathsf T}=(10/9)I_9$.
Thus
\[
  \sum_{i=1}^{10}(\ell_i^{\mathsf T}A\ell_i)(r_i^{\mathsf T}Br_i)
  =\tr\bigl((A\otimes B)ZZ^{\mathsf T}\bigr)
  =\frac{10}{9}\tr(A)\tr(B).
\]
Taking $B=I_3$ with arbitrary $A$, and then $A=I_3$ with arbitrary $B$, gives \eqref{eq:pair-frame}.

Let $R$ be the Gram matrix of $(r_i)_{i=1}^{10}$, and let $e_1,\ldots,e_{10}$ be the standard basis of $\mathbb R^{10}$.
Suppose $\ell_i=\varepsilon\ell_j$ for some $i\neq j$ and $\varepsilon\in\{1,-1\}$.
Then \eqref{eq:pair-condition} gives $R_{ij}=-\varepsilon/9$ and $R_{ki}=\varepsilon R_{kj}$ for $k\notin\{i,j\}$, since \eqref{eq:pair-condition} gives $\langle\ell_k,\ell_i\rangle\neq0$.
Hence
\[
  R(e_i-\varepsilon e_j)=\frac{10}{9}(e_i-\varepsilon e_j).
\]
On the other hand, \eqref{eq:pair-frame} gives $R^2=(10/3)R$, so the only nonzero eigenvalue of $R$ is $10/3$, a contradiction.
Thus $\ell_i\neq\pm\ell_j$ for $i\neq j$.
Equation \eqref{eq:pair-condition} is unchanged when $\ell_i$ and $r_i$ are interchanged, so $r_i\neq\pm r_j$ as well.
\end{proof}

\begin{lemma}
\label[lemma]{lem:pair-conic}
Under the assumptions of \cref{lem:pair-moments}, after exchanging $\ell_i$ and $r_i$ for every $i$ if necessary, there exists an invertible $H\in\Sym(3)$ such that
\[
  \tr H=0,\qquad \ell_i^{\mathsf T}H\ell_i=0\quad(1\leq i\leq10).
\]
\end{lemma}

\begin{proof}
By \eqref{eq:pair-frame} and \eqref{eq:pair-quadratic}, the subspaces
\[
  \{(\ell_i^{\mathsf T}A\ell_i)_{i=1}^{10}:A\in\Sym(3),\ \tr A=0\},
  \qquad
  \{(r_i^{\mathsf T}Br_i)_{i=1}^{10}:B\in\Sym(3),\ \tr B=0\}
\]
of $\mathbb R^{10}$ lie in $(1,\ldots,1)^\perp$ and are orthogonal to one another.
Their dimensions therefore sum to at most $9$.
Since the space of traceless matrices in $\Sym(3)$ has dimension $5$, after interchanging $\ell_i$ and $r_i$ for every $i$ if necessary, we can choose $H\in\Sym(3)\setminus\{0\}$ with $\tr H=0$ and $\ell_i^{\mathsf T}H\ell_i=0$ for every $i$.

A nonzero traceless symmetric $3\times3$ matrix has rank $2$ or $3$.
Suppose $H$ has rank $2$.
Its two nonzero eigenvalues are opposite, so rotating their orthonormal eigenvectors through $\pi/4$ gives coordinates in which $\ell_i^{\mathsf T}H\ell_i=0$ is $\ell_{i,1}\ell_{i,2}=0$.
Here $\ell_i=(\ell_{i,1},\ell_{i,2},\ell_{i,3})^{\mathsf T}$.
Put
\[
  I:=\{i\in\{1,\ldots,10\}:\ell_{i,2}=0\},\qquad
  J:=\{1,\ldots,10\}\setminus I.
\]
Then $\ell_{j,1}=0$ for $j\in J$.
Taking $A=\operatorname{diag}(1,0,0)$ and $A=\operatorname{diag}(0,1,0)$ in \eqref{eq:pair-quadratic}, with arbitrary $B\in\Sym(3)$, gives
\[
  \sum_{i\in I}\ell_{i,1}^2r_ir_i^{\mathsf T}
  =\sum_{j\in J}\ell_{j,2}^2r_jr_j^{\mathsf T}
  =\frac{10}{9}I_3.
\]
Thus each of $\{r_i:i\in I\}$ and $\{r_j:j\in J\}$ spans $\mathbb R^3$.
For $i\in I$ and $j\in J$, however, \eqref{eq:pair-condition} gives
\[
  \langle\ell_i,\ell_j\rangle=\ell_{i,3}\ell_{j,3},\qquad
  \langle r_i,r_j\rangle=-\frac1{9\ell_{i,3}\ell_{j,3}}.
\]
The product $\ell_{i,3}\ell_{j,3}$ is nonzero because $i\neq j$ and the product of the two inner products is $-1/9$.
Therefore the matrix $(\langle r_i,r_j\rangle)_{i\in I,j\in J}$ has rank $1$.
If $U\in\mathbb R^{3\times|I|}$ and $V\in\mathbb R^{3\times|J|}$ have columns $r_i$ $(i\in I)$ and $r_j$ $(j\in J)$, respectively, then $V$ is surjective and $U^{\mathsf T}$ is injective, so $U^{\mathsf T}V$ has rank $3$.
This contradiction shows that $H$ is invertible.
\end{proof}

\begin{lemma}
\label[lemma]{lem:finite-conic}
Under the assumptions of \cref{lem:pair-moments}, after possibly interchanging $\ell_i$ and $r_i$ for every $i$, applying a common orthogonal transformation to the $\ell_i$, and replacing individual pairs $(\ell_i,r_i)$ by $(-\ell_i,-r_i)$, there exist $b>0$ and distinct real numbers $t_1,\ldots,t_{10}$ such that
\[
  w(t):=\left(1,t,\frac{t^2-1}{b}\right)^{\mathsf T}\quad(t\in\mathbb R),
  \qquad
  \ell_i=\frac{w(t_i)}{\lVert w(t_i)\rVert}\quad(1\leq i\leq10).
\]
\end{lemma}

\begin{proof}
By \cref{lem:pair-conic}, it suffices to treat the case in which an invertible $H\in\Sym(3)$ satisfies $\tr H=0$ and $\ell_i^{\mathsf T}H\ell_i=0$ for every $i$.
One of $H$ and $-H$ has two positive eigenvalues, so choose $H$ with this property and denote its positive eigenvalues by $\lambda_1\geq\lambda_2>0$.
The third eigenvalue is $-\lambda_1-\lambda_2$ because $\tr H=0$.
There exists $U\in O(3)$ such that
\[
  U^{\mathsf T}HU
  =\operatorname{diag}(\lambda_1,\lambda_2,-\lambda_1-\lambda_2).
\]
Put $m:=1+2\lambda_1/\lambda_2\geq3$ and define
\[
  H_0:=\frac2{\lambda_2}U^{\mathsf T}HU
  =\operatorname{diag}(m-1,2,-m-1).
\]
Then $(U^{\mathsf T}\ell_i)^{\mathsf T}H_0(U^{\mathsf T}\ell_i)=0$ for every $i$.
If the third coordinate of $U^{\mathsf T}\ell_i$ were zero, the positivity of the first two diagonal entries of $H_0$ would force its other coordinates to be zero, contrary to $\lVert\ell_i\rVert=1$.
Thus there is a unique $\sigma_i\in\{1,-1\}$ such that the vectors
\[
  \widehat\ell_i:=\sigma_iU^{\mathsf T}\ell_i,
  \qquad
  \widehat r_i:=\sigma_i r_i
\]
satisfy $\widehat\ell_{i,3}>0$.
They are unit vectors and satisfy $\widehat\ell_i^{\mathsf T}H_0\widehat\ell_i=0$.
These operations preserve \eqref{eq:pair-condition}, so \cref{lem:pair-moments} also applies to $(\widehat\ell_i,\widehat r_i)$.

Consider the finite set
\[
  \mathcal C:=\{x=(x_1,x_2,x_3)^{\mathsf T}\in S^2:
  x^{\mathsf T}H_0x=0,\ x_1x_2=0,\ x_3>0\}.
\]
\[
  \mathcal C=\{a_1,a_{-1},d_1,d_{-1}\},
\]
where, for $\varepsilon\in\{1,-1\}$,
\[
  a_\varepsilon:=
  \left(\varepsilon\sqrt{\frac{m+1}{2m}},0,
  \sqrt{\frac{m-1}{2m}}\right)^{\mathsf T},
  \qquad
  d_\varepsilon:=
  \left(0,\varepsilon\sqrt{\frac{m+1}{m+3}},
  \sqrt{\frac2{m+3}}\right)^{\mathsf T}.
\]
We show that at least one point of $\mathcal C$ is not among the $\widehat\ell_i$.
Otherwise, there would be distinct indices $i_1,i_2,i_3,i_4$ with
\[
  (\widehat\ell_{i_1},\widehat\ell_{i_2},
   \widehat\ell_{i_3},\widehat\ell_{i_4})
  =(a_1,a_{-1},d_1,d_{-1}).
\]
The inner products among $a_1,a_{-1},d_1,d_{-1}$ are
\[
  \begin{aligned}
  \langle a_1,a_{-1}\rangle&=-\frac1m,
  \qquad \langle d_1,d_{-1}\rangle=-\frac{m-1}{m+3},\\
  \langle a_\varepsilon,d_\delta\rangle
  &=\sqrt{\frac{m-1}{m(m+3)}}
  \quad(\varepsilon,\delta\in\{1,-1\}).
  \end{aligned}
\]
Set
\[
  \alpha:=\frac m9,
  \qquad
  \beta:=\frac{m+3}{9(m-1)},
  \qquad
  z:=(\sqrt\alpha,\sqrt\alpha,-\sqrt\beta,-\sqrt\beta)^{\mathsf T}.
\]
Equation \eqref{eq:pair-condition} now shows that the Gram matrix of $(\widehat r_{i_s})_{s=1}^4$ is
\[
  G=\operatorname{diag}(1-\alpha,1-\alpha,1-\beta,1-\beta)
    +zz^{\mathsf T}.
\]
Here $\alpha,\beta>0$ because $m\geq3$.
By \cref{lem:pair-moments}, the unit vectors $\widehat r_{i_1},\widehat r_{i_2}$ are distinct, as are $\widehat r_{i_3},\widehat r_{i_4}$, so
\[
  1-\alpha
  =\frac12\lVert\widehat r_{i_1}-\widehat r_{i_2}\rVert^2>0,
  \qquad
  1-\beta
  =\frac12\lVert\widehat r_{i_3}-\widehat r_{i_4}\rVert^2>0.
\]
Since $1-\alpha,1-\beta>0$ and $zz^{\mathsf T}$ is positive semidefinite, $G$ is positive definite and has rank $4$, contrary to its being the Gram matrix of four vectors in $\mathbb R^3$.

Choose $p\in\mathcal C\setminus\{\widehat\ell_i:1\leq i\leq10\}$.
Also $\widehat\ell_i\neq-p$ for every $i$, since $p$ and every $\widehat\ell_i$ have positive third coordinate.
Put $v:=(0,1,0)^{\mathsf T}$ if $p_2=0$, and $v:=(1,0,0)^{\mathsf T}$ otherwise.
Then $p,v$ are orthogonal unit vectors and $H_0v=\lambda v$ for some $\lambda\in\{2,m-1\}$, with $\lambda>0$.
Choose a unit vector $q$ perpendicular to both $p$ and $v$, and put $V:=(q\ v\ p)\in O(3)$.
\[
  q^{\mathsf T}H_0v=0,
  \qquad
  v^{\mathsf T}H_0p=0,
  \qquad
  v^{\mathsf T}H_0v=\lambda.
\]
Since $p^{\mathsf T}H_0p=0$ and $\tr H_0=0$, it follows that
\[
  V^{\mathsf T}H_0V
  =\begin{pmatrix}
    -\lambda&0&c\\
    0&\lambda&0\\
    c&0&0
  \end{pmatrix},
  \qquad c:=q^{\mathsf T}H_0p.
\]
Since $H_0$ is invertible, $c\neq0$.
Replacing $q$ by $-q$ if necessary, we may assume $c<0$, and put $b:=-2c/\lambda>0$.
Then
\[
  V^{\mathsf T}H_0V
  =\frac\lambda2
  \begin{pmatrix}
    -2&0&-b\\
    0&2&0\\
    -b&0&0
  \end{pmatrix}.
\]

For each $i$, put $u_i:=V^{\mathsf T}\widehat\ell_i=(u_{i,1},u_{i,2},u_{i,3})^{\mathsf T}$.
The vector $u_i$ has norm $1$, and
\[
  0=u_i^{\mathsf T}(V^{\mathsf T}H_0V)u_i
  =\lambda(-u_{i,1}^2+u_{i,2}^2-bu_{i,1}u_{i,3}).
\]
Thus $u_{i,2}^2=u_{i,1}^2+bu_{i,1}u_{i,3}$.
If $u_{i,1}=0$, then $u_{i,2}=0$ and $\lVert u_i\rVert=1$ gives $u_i=\pm(0,0,1)^{\mathsf T}$.
Since the third column of $V$ is $p$, this would imply $\widehat\ell_i=Vu_i=\pm p$, contrary to the choice of $p$.
Therefore $u_{i,1}\neq0$, and we can define $t_i:=u_{i,2}/u_{i,1}\in\mathbb R$.
\[
  \frac{u_{i,3}}{u_{i,1}}=\frac{t_i^2-1}{b},
  \qquad
  u_i=u_{i,1}\left(1,t_i,\frac{t_i^2-1}{b}\right)^{\mathsf T}
     =u_{i,1}w(t_i).
\]
Choose $\varepsilon_i\in\{1,-1\}$ so that $\varepsilon_i u_{i,1}>0$.
Then $\varepsilon_i u_{i,1}=|u_{i,1}|$, and $\lVert u_i\rVert=1$ gives $|u_{i,1}|\lVert w(t_i)\rVert=1$.
Define
\[
  \widetilde\ell_i:=\varepsilon_iV^{\mathsf T}\widehat\ell_i,
  \qquad
  \widetilde r_i:=\varepsilon_i\widehat r_i.
\]
Then
\[
  \widetilde\ell_i
  =\varepsilon_i u_i
  =|u_{i,1}|w(t_i)
  =\frac{w(t_i)}{\lVert w(t_i)\rVert}.
\]
If $t_i=t_j$ for distinct $i,j$, then $\widetilde\ell_i=\widetilde\ell_j$, and hence $\widehat\ell_i=\varepsilon_i\varepsilon_j\widehat\ell_j$, contrary to \cref{lem:pair-moments}.
Thus $t_1,\ldots,t_{10}$ are distinct.

\end{proof}

We will apply the following lemma to obtain $b^2=8$ and a formula for $\prod_{i=1}^{10}(t-t_i)$.
Here $\mathrm i$ denotes the imaginary unit.

\begin{lemma}
\label[lemma]{lem:multiplication-matrix}
Let $b>0$, $K=K^{\mathsf T}\in\mathbb R^{3\times3}$, $u\in\mathbb R^3$, and $a\in\mathbb R$, and put
\begin{equation}
  M:=\begin{pmatrix}
    0&I_3&0&0\\
    I_3&0&bI_3&0\\
    0&bI_3&K&u\\
    0&0&u^{\mathsf T}&a
  \end{pmatrix}.
  \label{eq:multiplication-matrix}
\end{equation}
Put $g(t):=\det(tI_{10}-M)$, and let $F(t)$ be the $(10,10)$-cofactor of $tI_{10}-M$.
If $g'=10F$, then
\[
  b^2=8,\qquad
  g(t)=\operatorname{Re}\bigl((t+\mathrm i\sqrt3)^{10}\bigr)
       -\frac a{\sqrt3}\operatorname{Im}\bigl((t+\mathrm i\sqrt3)^{10}\bigr).
\]
\end{lemma}

\begin{proof}
We first compute $g$ and $F$ using an orthogonal diagonalization of $K$.
Choose $U\in O(3)$ such that
\[
  U^{\mathsf T}KU=\Lambda,
  \qquad \Lambda:=\operatorname{diag}(\lambda_1,\lambda_2,\lambda_3).
\]
Put $\widetilde u:=U^{\mathsf T}u=(\eta_1,\eta_2,\eta_3)^{\mathsf T}$, $S_9:=\operatorname{diag}(U,U,U)$, and $S:=\operatorname{diag}(S_9,1)$.
Then
\[
  \widetilde M:=S^{\mathsf T}MS
  =\begin{pmatrix}
    0&I_3&0&0\\
    I_3&0&bI_3&0\\
    0&bI_3&\Lambda&\widetilde u\\
    0&0&\widetilde u^{\mathsf T}&a
  \end{pmatrix}.
\]
Thus $\widetilde M$ has the block form \eqref{eq:multiplication-matrix}, since $U^{\mathsf T}I_3U=I_3$ and $U^{\mathsf T}(bI_3)U=bI_3$.
We next check that conjugation by $S$ preserves $g$ and $F$.
Orthogonality of $S$ gives
\[
  \det(tI_{10}-\widetilde M)
  =\det\bigl(S^{\mathsf T}(tI_{10}-M)S\bigr)
  =(\det S)^2g(t)=g(t).
\]
Let $B:=(M_{ij})_{1\leq i,j\leq9}$ and $\widetilde B:=(\widetilde M_{ij})_{1\leq i,j\leq9}$.
Since the sign of the $(10,10)$-cofactor is $(-1)^{20}=1$, we have $F(t)=\det(tI_9-B)$.
Moreover, $\widetilde B=S_9^{\mathsf T}BS_9$, so
\[
  \det(tI_9-\widetilde B)
  =\det\bigl(S_9^{\mathsf T}(tI_9-B)S_9\bigr)
  =(\det S_9)^2F(t)=F(t).
\]
We next introduce three $3\times3$ matrices to compute $F$ and $g$.
Define
\[
  k:=1+b^2>1,\qquad P(t):=t^3-kt,\qquad Q(t):=t^2-1.
\]
For $1\leq j\leq3$, put
\[
  A_j(t):=\begin{pmatrix}
    t&-1&0\\
    -1&t&-b\\
    0&-b&t-\lambda_j
  \end{pmatrix},
  \qquad d_j(t):=\det A_j(t).
\]
Expansion along the first row gives
\[
  \begin{aligned}
  d_j(t)
  &=t\bigl(t(t-\lambda_j)-b^2\bigr)-(t-\lambda_j)\\
  &=t^3-(1+b^2)t-\lambda_j(t^2-1)\\
  &=P(t)-\lambda_jQ(t).
  \end{aligned}
\]
We first compute $F$ by grouping the rows and columns into the blocks $A_1(t),A_2(t),A_3(t)$.
Let $e_3:=(0,0,1)^{\mathsf T}$.
Simultaneously reordering the rows and columns of $tI_{10}-\widetilde M$ in the order $(1,4,7,2,5,8,3,6,9,10)$ gives
\[
  \begin{pmatrix}
    A_1(t)&0&0&-\eta_1e_3\\
    0&A_2(t)&0&-\eta_2e_3\\
    0&0&A_3(t)&-\eta_3e_3\\
    -\eta_1e_3^{\mathsf T}&-\eta_2e_3^{\mathsf T}&-\eta_3e_3^{\mathsf T}&t-a
  \end{pmatrix}.
\]
The rows and columns are permuted together, so the determinant is unchanged.
As the last row and column remain last, the determinant of the leading $9\times9$ submatrix is also unchanged, and hence
\begin{equation}
  F(t)=\prod_{j=1}^3d_j(t)=\prod_{j=1}^3\bigl(P(t)-\lambda_jQ(t)\bigr).
  \label{eq:multiplication-F-product}
\end{equation}

We next compute $g$.
For $d_j(t)\neq0$, the cofactor formula gives
\[
  e_3^{\mathsf T}A_j(t)^{-1}e_3
  =\frac{\det\begin{pmatrix}t&-1\\-1&t\end{pmatrix}}{d_j(t)}
  =\frac{Q(t)}{d_j(t)}.
\]
Put
\[
  \mathcal B(t):=\operatorname{diag}(A_1(t),A_2(t),A_3(t)),
  \qquad
  \xi:=\begin{pmatrix}-\eta_1e_3\\-\eta_2e_3\\-\eta_3e_3\end{pmatrix}.
\]
Subtracting $\xi^{\mathsf T}\mathcal B(t)^{-1}$ times the first nine rows from the last row gives
\[
  \begin{aligned}
  g(t)
  &=\det\begin{pmatrix}\mathcal B(t)&\xi\\\xi^{\mathsf T}&t-a\end{pmatrix}\\
  &=\det\mathcal B(t)\bigl(t-a-\xi^{\mathsf T}\mathcal B(t)^{-1}\xi\bigr)\\
  &=F(t)\left(t-a-\sum_{j=1}^3\eta_j^2e_3^{\mathsf T}A_j(t)^{-1}e_3\right)\\
  &=F(t)\left(t-a-Q(t)\sum_{j=1}^3\frac{\eta_j^2}{d_j(t)}\right).
  \end{aligned}
\]
Consequently,
\begin{equation}
  g(t)=(t-a)F(t)-Q(t)\sum_{j=1}^3\eta_j^2
  \prod_{\substack{1\leq h\leq3\\h\neq j}}
  \bigl(P(t)-\lambda_hQ(t)\bigr).
  \label{eq:multiplication-g-formula}
\end{equation}
Each $d_j$ is a monic cubic polynomial, so $d_1(t)d_2(t)d_3(t)=0$ for only finitely many real $t$.
Both sides of \eqref{eq:multiplication-g-formula} are polynomials and agree whenever $d_1(t)d_2(t)d_3(t)\neq0$, so \eqref{eq:multiplication-g-formula} holds for every $t$.

We next use the even part of $g$ to prove that $b^2=8$.
Separate $g$ into its even and odd parts by defining
\[
  g_{\mathrm{ev}}(t):=\frac{g(t)+g(-t)}2,
  \qquad
  g_{\mathrm{odd}}(t):=\frac{g(t)-g(-t)}2.
\]
Also put
\[
  \begin{aligned}
  \beta&:=\sum_{1\leq j<h\leq3}\lambda_j\lambda_h,\\
  \alpha&:=a\sum_{j=1}^3\lambda_j-\sum_{j=1}^3\eta_j^2,\\
  \gamma&:=a\lambda_1\lambda_2\lambda_3
  -\sum_{j=1}^3\eta_j^2
    \prod_{\substack{1\leq h\leq3\\h\neq j}}\lambda_h.
  \end{aligned}
\]
We first express $g_{\mathrm{ev}}$ in terms of $P$ and $Q$.
We omit the argument $t$ in $P,Q,F,g_{\mathrm{ev}}$, and their derivatives.
Expanding the product for $F$ gives
\[
  F=P^3-\left(\sum_{j=1}^3\lambda_j\right)P^2Q
    +\beta PQ^2-\lambda_1\lambda_2\lambda_3Q^3.
\]
Since $P$ is odd and $Q$ is even, the even part of $(t-a)F$ is
\[
  tP^3+\beta tPQ^2
  +a\left(\sum_{j=1}^3\lambda_j\right)P^2Q
  +a\lambda_1\lambda_2\lambda_3Q^3.
\]
For each $j$, expansion gives
\[
  Q\prod_{\substack{1\leq h\leq3\\h\neq j}}(P-\lambda_hQ)
  =P^2Q-\left(\sum_{\substack{1\leq h\leq3\\h\neq j}}\lambda_h\right)PQ^2
   +\left(\prod_{\substack{1\leq h\leq3\\h\neq j}}\lambda_h\right)Q^3.
\]
The even part of $Q\prod_{h\neq j}(P-\lambda_hQ)$ is $P^2Q+\bigl(\prod_{h\neq j}\lambda_h\bigr)Q^3$.
Thus the definitions of $\alpha,\gamma$ and \eqref{eq:multiplication-g-formula} give
\begin{equation}
  g_{\mathrm{ev}}=tP^3+\beta tPQ^2+\alpha P^2Q+\gamma Q^3.
  \label{eq:multiplication-even-polynomial}
\end{equation}
We next derive two expressions for $g_{\mathrm{ev}}'$ and compare them to determine $\alpha,\beta,\gamma$.
Taking the odd part of $g'=10F$, we obtain
\begin{equation}
  g_{\mathrm{ev}}'=10(P^3+\beta PQ^2).
  \label{eq:multiplication-even-derivative}
\end{equation}
Differentiating \eqref{eq:multiplication-even-polynomial} gives
\begin{equation}
  \begin{aligned}
  g_{\mathrm{ev}}'
  ={}&P^3+3tP^2P'
  +\beta\bigl(PQ^2+tP'Q^2+2tPQQ'\bigr)\\
  &+\alpha\bigl(2PP'Q+P^2Q'\bigr)+3\gamma Q^2Q'.
  \end{aligned}
  \label{eq:multiplication-even-derivative-expanded}
\end{equation}
To determine $\alpha$, first set $t=1$.
At this point, $Q(1)=0$, $Q'(1)=2$, $P(1)=1-k\neq0$, and $P'(1)=3-k$.
Evaluating \cref{eq:multiplication-even-derivative,eq:multiplication-even-derivative-expanded} at $t=1$ and dividing by $P(1)^2$ gives
\[
  3(3-k)+2\alpha=9(1-k),
  \qquad\text{hence}\qquad \alpha=-3k.
\]
Next, we determine $\beta$ by comparing the coefficients of $t^7$.
The coefficients of $t^8$ in $g_{\mathrm{ev}}$ and of $t^7$ in $P^3+\beta PQ^2$ are $-3k+\beta+\alpha$ and $-3k+\beta$, respectively.
Comparing the coefficients of $t^7$ in \eqref{eq:multiplication-even-derivative} gives
\[
  8(-3k+\beta+\alpha)=10(-3k+\beta).
\]
Using $\alpha=-3k$, we obtain $\beta=-9k$.
Next, we determine $\gamma$ by evaluating at a zero of $P$.
Put $t_0:=\sqrt{k}$, so $P(t_0)=0$, $P'(t_0)=2k$, $Q(t_0)=k-1$, and $Q'(t_0)=2t_0$.
Evaluating \cref{eq:multiplication-even-derivative,eq:multiplication-even-derivative-expanded} at $t=t_0$ gives
\[
  \begin{aligned}
  0=g_{\mathrm{ev}}'(t_0)
  &=\beta t_0P'(t_0)Q(t_0)^2+3\gamma Q(t_0)^2Q'(t_0)\\
  &=2t_0(k-1)^2(k\beta+3\gamma).
  \end{aligned}
\]
Since $t_0>0$ and $k>1$, we obtain $k\beta+3\gamma=0$, and hence $\gamma=3k^2$.

Substituting $\alpha=-3k$, $\beta=-9k$, and $\gamma=3k^2$ into \eqref{eq:multiplication-even-polynomial} gives
\[
  g_{\mathrm{ev}}=tP^3-9ktPQ^2-3kP^2Q+3k^2Q^3.
\]
Using \eqref{eq:multiplication-even-derivative}, $Q'=2t$, and $tP'-kQ'=3P$, we obtain
\[
  \begin{aligned}
  0
  &=g_{\mathrm{ev}}'-10(P^3-9kPQ^2)\\
  &=3P^2(tP'-3P-kQ')+81kPQ^2+9kQ^2(kQ'-tP')\\
  &\qquad{}-18ktPQQ'-6kPP'Q\\
  &=54kPQ^2-36kt^2PQ-6kPP'Q\\
  &=6kPQ(9Q-6t^2-P')\\
  &=6k(k-9)PQ.
  \end{aligned}
\]
Here $9Q-6t^2-P'=9(t^2-1)-6t^2-(3t^2-k)=k-9$.
Since $P,Q$ are nonzero polynomials and $k>1$, we conclude that $k=9$, and hence $b^2=8$.

We next identify $g_{\mathrm{ev}}$ with $\operatorname{Re}((t+\mathrm i\sqrt3)^{10})$.
With $k=9$, direct expansion gives
\begin{equation}
  (t+\mathrm i\sqrt3)^3
  =(t^3-9t)+3\mathrm i\sqrt3(t^2-1)
  =P+3\mathrm i\sqrt3Q
  \label{eq:multiplication-cubic-expansion}
\end{equation}
and
\[
  (P+3\mathrm i\sqrt3Q)^3
     =P^3-81PQ^2+\mathrm i\sqrt3(9P^2Q-81Q^3).
\]
Therefore
\[
  \begin{aligned}
  g_{\mathrm{ev}}
  &=t(P^3-81PQ^2)-3(9P^2Q-81Q^3)\\
  &=\operatorname{Re}\bigl((t+\mathrm i\sqrt3)(P+3\mathrm i\sqrt3Q)^3\bigr)\\
  &=\operatorname{Re}\bigl((t+\mathrm i\sqrt3)^{10}\bigr).
  \end{aligned}
\]

We next determine the odd part of $g$.
We first determine $\lambda_1+\lambda_2+\lambda_3$ using \cref{eq:multiplication-F-product,eq:multiplication-g-formula}.
Each factor $P(t)-\lambda_jQ(t)=t^3-\lambda_jt^2-kt+\lambda_j$ is a monic cubic, so $F$ is monic of degree $9$ and its coefficient of $t^8$ is $-\sum_{j=1}^3\lambda_j$.
Since $Q$ has degree $2$ and each product $\prod_{h\neq j}(P-\lambda_hQ)$ has degree $6$, the term subtracted from $(t-a)F$ in \eqref{eq:multiplication-g-formula} has degree at most $8$.
Hence the coefficient of $t^9$ in $g$ is $-\sum_{j=1}^3\lambda_j-a$.
Comparing the coefficients of $t^8$ in $g'=10F$ now gives
\[
  -9\left(\sum_{j=1}^3\lambda_j+a\right)
  =-10\sum_{j=1}^3\lambda_j,
  \qquad\text{hence}\qquad
  \sum_{j=1}^3\lambda_j=9a.
\]
Next, we compute $g_{\mathrm{odd}}'$ using the even part of $F$.
Using $\lambda_1+\lambda_2+\lambda_3=9a$ and differentiating $g_{\mathrm{odd}}(t)=(g(t)-g(-t))/2$, we obtain
\[
  \begin{aligned}
  g_{\mathrm{odd}}'(t)
  &=\frac{g'(t)+g'(-t)}2=5\bigl(F(t)+F(-t)\bigr)\\
  &=-10\left(\sum_{j=1}^3\lambda_j\right)P^2Q-10\lambda_1\lambda_2\lambda_3Q^3\\
  &=-90aP^2Q-10\lambda_1\lambda_2\lambda_3Q^3.
  \end{aligned}
\]
We next prove that $g_{\mathrm{odd}}$ equals $-(a/\sqrt3)\operatorname{Im}((t+\mathrm i\sqrt3)^{10})$.
To do this, define their difference by
\[
  h(t):=g_{\mathrm{odd}}(t)
       +\frac a{\sqrt3}\operatorname{Im}\bigl((t+\mathrm i\sqrt3)^{10}\bigr).
\]
We first compute $h'$.
Equation \eqref{eq:multiplication-cubic-expansion} gives
\[
  \frac{d}{dt}\operatorname{Im}\bigl((t+\mathrm i\sqrt3)^{10}\bigr)
  =10\operatorname{Im}\bigl((t+\mathrm i\sqrt3)^9\bigr)
  =90\sqrt3(P^2Q-9Q^3).
\]
Consequently,
\[
  \begin{aligned}
  h'
  &=-90aP^2Q-10\lambda_1\lambda_2\lambda_3Q^3
     +90a(P^2Q-9Q^3)\\
  &=\mu Q^3,
  \qquad\text{where }\mu:=-10\lambda_1\lambda_2\lambda_3-810a.
  \end{aligned}
\]
We next show that $h(0)=h(1)=0$.
Since $g_{\mathrm{odd}}(0)=0$ and $(\mathrm i\sqrt3)^{10}=-3^5$ is real, we have $h(0)=0$.
To evaluate $h(1)$, use $Q(\pm1)=0$ in \eqref{eq:multiplication-g-formula} to obtain
\[
  g(\pm1)=(\pm1-a)P(\pm1)^3.
\]
As $P(-1)=-P(1)$, it follows that
\[
  g_{\mathrm{odd}}(1)
  =\frac{(1-a)P(1)^3-(1+a)P(1)^3}{2}
  =-aP(1)^3.
\]
Also, $(1+\mathrm i\sqrt3)^3=P(1)$, and hence
\[
  (1+\mathrm i\sqrt3)^{10}
  =(1+\mathrm i\sqrt3)P(1)^3,
  \qquad
  \operatorname{Im}\bigl((1+\mathrm i\sqrt3)^{10}\bigr)=\sqrt3P(1)^3.
\]
Hence $h(1)=-aP(1)^3+(a/\sqrt3)\sqrt3P(1)^3=0$.
Finally, we use $h(0)=h(1)=0$ to show that $h=0$.
Rolle's theorem applied to $h$ on $[0,1]$ gives some $t_*\in(0,1)$ with $h'(t_*)=0$.
Since $Q(t_*)=t_*^2-1<0$, the equality $0=\mu Q(t_*)^3$ implies $\mu=0$.
Therefore $h'=0$, and $h(0)=0$ gives $h=0$.
Combining the even and odd parts yields
\[
  g(t)=\operatorname{Re}\bigl((t+\mathrm i\sqrt3)^{10}\bigr)
       -\frac a{\sqrt3}\operatorname{Im}\bigl((t+\mathrm i\sqrt3)^{10}\bigr),
\]
as required.
\end{proof}

\begin{proof}[Proof of \cref{thm:tight-uniqueness}]
Let $D\subset G_{2,4}$ be a tight Grassmannian $4$-design, and write
\[
  D=\{[\ell_i,r_i]:1\leq i\leq10\},
  \qquad \ell_i,r_i\in S^2.
\]
The design satisfies \eqref{eq:pair-condition}.
Since quaternionic conjugation belongs to $O(4)$ and sends $[\ell,r]$ to $[r,\ell]$, we may interchange $\ell_i$ and $r_i$ for every $i$ without loss of generality.
Let $Z\in\mathbb R^{9\times10}$ have columns $\ell_i\otimes r_i$.
By \cref{lem:G24-moment-criterion},
\begin{equation}
  Z(1,\ldots,1)^{\mathsf T}=0,
  \qquad
  ZZ^{\mathsf T}=\frac{10}{9}I_9.
  \label{eq:uniqueness-matrix-identities}
\end{equation}

\Cref{lem:finite-conic} gives $U\in O(3)$, signs $\sigma_i\in\{1,-1\}$, a number $b>0$, and distinct real numbers $t_1,\ldots,t_{10}$ such that the vectors
\[
  \widehat\ell_i:=\sigma_iU^{\mathsf T}\ell_i,
  \qquad
  \widehat r_i:=\sigma_i r_i
\]
satisfy
\[
  \widehat\ell_i=\frac{w(t_i)}{\lVert w(t_i)\rVert}.
\]

Let $\widehat Z$ have columns $\widehat\ell_i\otimes\widehat r_i$.
Since $\sigma_i^2=1$, we have
\[
  \widehat Z=(U^{\mathsf T}\otimes I_3)Z.
\]
Since $U^{\mathsf T}\otimes I_3$ is orthogonal, \eqref{eq:uniqueness-matrix-identities} gives
\begin{equation}
  \widehat Z(1,\ldots,1)^{\mathsf T}=0,
  \qquad
  \widehat Z\widehat Z^{\mathsf T}=\frac{10}{9}I_9.
  \label{eq:uniqueness-transformed-moments}
\end{equation}
The transformed pairs also satisfy
\[
  \langle\widehat\ell_i,\widehat\ell_j\rangle
  \langle\widehat r_i,\widehat r_j\rangle
  =\frac{10\delta_{ij}-1}{9}.
\]
Indeed, $U$ preserves inner products, and the factor $\sigma_i\sigma_j$ occurs in both inner products, so its contribution to their product is $1$.

To apply \cref{lem:multiplication-matrix}, put
\[
  T:=\operatorname{diag}(t_1,\ldots,t_{10}),
\]
and write $\widehat\ell_{i,k},\widehat r_{i,j}$ for the coordinates of the transformed vectors.
Define
\[
  e_{kj}:=\sqrt{\frac9{10}}
  \bigl(\widehat\ell_{i,k}\widehat r_{i,j}\bigr)_{i=1}^{10}
  \quad(1\leq k,j\leq3),
  \qquad
  e_*:=\frac1{\sqrt{10}}(1,\ldots,1)^{\mathsf T}.
\]
Equation \eqref{eq:uniqueness-transformed-moments} gives
\[
  \begin{aligned}
  \langle e_{kj},e_{mn}\rangle
  &=\frac9{10}\sum_{i=1}^{10}
    \widehat\ell_{i,k}\widehat\ell_{i,m}
    \widehat r_{i,j}\widehat r_{i,n}
   =\delta_{km}\delta_{jn},\\
  \langle e_{kj},e_*\rangle
  &=\frac3{10}\sum_{i=1}^{10}
    \widehat\ell_{i,k}\widehat r_{i,j}
   =0.
  \end{aligned}
\]
Since $\lVert e_*\rVert=1$, the vectors $e_{kj}$ $(1\leq k,j\leq3)$ and $e_*$ form an orthonormal basis of $\mathbb R^{10}$.

The coordinates of $\widehat\ell_i=w(t_i)/\lVert w(t_i)\rVert$ satisfy
\[
  \widehat\ell_{i,2}=t_i\widehat\ell_{i,1},
  \qquad
  t_i\widehat\ell_{i,2}
  =\widehat\ell_{i,1}+b\widehat\ell_{i,3}.
\]
Multiplying the coordinates in the definitions of $e_{1j}$ and $e_{2j}$ by $t_i$ therefore gives
\[
  Te_{1j}=e_{2j},
  \qquad
  Te_{2j}=e_{1j}+be_{3j}
  \qquad(1\leq j\leq3).
\]

Order the basis as
\[
  e_{11},e_{12},e_{13},\;
  e_{21},e_{22},e_{23},\;
  e_{31},e_{32},e_{33},\;
  e_*,
\]
and let $M$ be the matrix of $T$ in this basis.
The relations $Te_{1j}=e_{2j}$ and $Te_{2j}=e_{1j}+be_{3j}$ determine the first two block columns of $M$.
Since $T$ is symmetric and the basis is orthonormal, $M$ is symmetric, so the corresponding block rows are also determined.
Thus
\[
  M=
  \begin{pmatrix}
    0&I_3&0&0\\
    I_3&0&bI_3&0\\
    0&bI_3&K&u\\
    0&0&u^{\mathsf T}&a
  \end{pmatrix}
\]
for some
\[
  K=K^{\mathsf T}\in\mathbb R^{3\times3},
  \qquad
  u\in\mathbb R^3,
  \qquad
  a\in\mathbb R.
\]

Put
\[
  g(t):=\det(tI_{10}-M),
\]
and let $F(t)$ be the determinant of the matrix obtained by deleting row $10$ and column $10$ from $tI_{10}-M$.
Since $M$ represents $T$ in an orthonormal basis,
\[
  g(t)=\det(tI_{10}-T)
  =\prod_{i=1}^{10}(t-t_i).
\]

For $t\notin\{t_1,\ldots,t_{10}\}$, we have
\[
  \frac{F(t)}{g(t)}
  =\bigl((tI_{10}-M)^{-1}\bigr)_{10,10}.
\]
The last vector of the chosen basis is $e_*$, so
\[
  \begin{aligned}
  \frac{F(t)}{g(t)}
  &=e_*^{\mathsf T}(tI_{10}-T)^{-1}e_*\\
  &=\frac1{10}\sum_{i=1}^{10}\frac1{t-t_i}\\
  &=\frac{g'(t)}{10g(t)}.
  \end{aligned}
\]
Hence $g'(t)=10F(t)$ for all $t$ outside a finite set.
Since both sides are polynomials,
\[
  g'=10F
\]
holds as a polynomial identity.

Applying \cref{lem:multiplication-matrix}, we obtain $b^2=8$ and
\[
  g(t)
  =\operatorname{Re}\bigl((t+\mathrm i\sqrt3)^{10}\bigr)
  -\frac a{\sqrt3}
   \operatorname{Im}\bigl((t+\mathrm i\sqrt3)^{10}\bigr).
\]

For each $i$, write
\[
  t_i=\sqrt3\cot\theta_i,
  \qquad
  \theta_i\in(0,\pi).
\]
Since
\[
  t_i+\mathrm i\sqrt3
  =\frac{\sqrt3}{\sin\theta_i}
   \bigl(\cos\theta_i+\mathrm i\sin\theta_i\bigr),
\]
$g(t_i)=0$ gives
\[
  0=\frac{3^5}{\sin^{10}\theta_i}
  \left(
    \cos(10\theta_i)
    -\frac a{\sqrt3}\sin(10\theta_i)
  \right).
\]
If $\sin(10\theta_i)=0$, then $\cos(10\theta_i)=\pm1$, contradicting $g(t_i)=0$.
Thus
\[
  \sin(10\theta_i)\neq0,
  \qquad
  \cot(10\theta_i)=\frac a{\sqrt3}.
\]

Let $\varphi\in(0,\pi)$ be the unique number satisfying $\cot\varphi=a/\sqrt3$.
Since $\cot$ has period $\pi$ and is injective on $(0,\pi)$, the solutions $\theta\in(0,\pi)$ of $\cot(10\theta)=a/\sqrt3$ form the set
\[
  \left\{
    \frac{\varphi+j\pi}{10}:0\leq j\leq9
  \right\}.
\]
The $t_i$ are distinct, so the $\theta_i$ are distinct.
After relabelling the pairs by $j\in\{0,\ldots,9\}$, we may assume that
\[
  \theta_j=\frac{\varphi+j\pi}{10}
  \qquad(0\leq j\leq9).
\]

Since $b^2=8$,
\[
  \lVert w(t)\rVert^2
  =1+t^2+\frac{(t^2-1)^2}{8}
  =\frac{(t^2+3)^2}{8}.
\]
Taking inner products of $w(t_j)/\lVert w(t_j)\rVert$ and substituting $t_j=\sqrt3\cot\theta_j$ and $\theta_j=(\varphi+j\pi)/10$, we obtain
\begin{equation}
  \begin{aligned}
  \langle\widehat\ell_j,\widehat\ell_k\rangle
  &=\frac{t_j^2t_k^2-t_j^2-t_k^2+8t_jt_k+9}{(t_j^2+3)(t_k^2+3)}\\
  &=\frac{1+2\cos\bigl(2(\theta_j-\theta_k)\bigr)}3\\
  &=\frac{1+2\cos((j-k)\pi/5)}3\\
  &=\langle u_j,u_k\rangle
  \end{aligned}
  \qquad(0\leq j,k\leq9).
  \label{eq:uniqueness-left-gram}
\end{equation}
By \eqref{eq:D10-Hadamard}, $\langle u_j,u_k\rangle\neq0$ for $j\neq k$.
Equations \eqref{eq:pair-condition} and \eqref{eq:uniqueness-left-gram} therefore give
\begin{equation}
  \begin{aligned}
  \langle\widehat r_j,\widehat r_k\rangle
  &=-\frac1{9\langle\widehat\ell_j,\widehat\ell_k\rangle}\\
  &=\langle v_j,v_k\rangle
  \end{aligned}
  \qquad(j\neq k).
  \label{eq:uniqueness-right-gram}
\end{equation}
For $j=k$, both $\langle\widehat r_j,\widehat r_k\rangle$ and $\langle v_j,v_k\rangle$ equal $1$.
Since $U\widehat\ell_j=\sigma_j\ell_j$ and $U$ is orthogonal, the family $(\sigma_j\ell_j)_{j=0}^9$ has the same Gram matrix as $(u_j)_{j=0}^9$.
Also, $\widehat r_j=\sigma_jr_j$, so $(\sigma_jr_j)_{j=0}^9$ has the same Gram matrix as $(v_j)_{j=0}^9$.
Thus the assignments $u_j\mapsto\sigma_j\ell_j$ and $v_j\mapsto\sigma_jr_j$ extend to well-defined linear maps between the respective spans, and both maps preserve norms.
By \cref{lem:pair-moments}, each of the families $(u_j)_j$, $(v_j)_j$, $(\ell_j)_j$, and $(r_j)_j$ spans $\mathbb R^3$.
Multiplying individual vectors by signs does not change their spans.
Hence these linear maps are represented by matrices $A,B\in O(3)$ satisfying
\[
  Au_j=\sigma_j\ell_j,
  \qquad
  Bv_j=\sigma_jr_j
  \qquad(0\leq j\leq9).
\]
Simultaneous sign changes leave the represented planes unchanged, so
\[
  D=\{[Au_j,Bv_j]:0\leq j\leq9\}.
\]

To represent $D$ using two matrices in $SO(3)$, put
\[
  C:=\operatorname{diag}(-1,-1,1).
\]
We have
\[
  u_{j+5}=Cu_j,
  \qquad
  v_{j+5}=-Cv_j.
\]
Shifting indices by $5$ therefore gives
\[
  D=\{[ACu_j,-BCv_j]:j\in\mathbb Z/10\mathbb Z\}.
\]
Also,
\[
  \det(AC)\det(-BC)=-\det A\det B.
\]
Define
\[
  (A_+,B_+):=
  \begin{cases}
    ((\det A)A,(\det B)B),&\text{if } \det(AB)=1,\\
    ((\det A)AC,(\det B)BC),&\text{if }\det(AB)=-1.
  \end{cases}
\]
Then $A_+,B_+\in SO(3)$ and
\[
  D=\{[A_+u_j,B_+v_j]:j\in\mathbb Z/10\mathbb Z\}.
\]
By \cref{sec:gram-criterion}, there is an element of $SO(4)$ sending $[\ell,r]$ to $[A_+\ell,B_+r]$ for all $\ell,r\in S^2$.
This element maps $D_{10}$ onto $D$, proving the theorem.
\end{proof}

\section{Primitive idempotents}
\label{sec:primitive-idempotents}

Let $(X,\{R_i\}_{i=0}^s)$ be a symmetric association scheme of order $N:=|X|$, with the notation of \cref{sec:association-scheme-construction}.
We now restrict the construction to primitive idempotents of rank $3$ and express the injectivity and design conditions in terms of the parameters of the scheme.
For $0\leq i\leq s$, the \emph{valency} of $R_i$ is
\[
  k_i:=|\{y\in X:(x,y)\in R_i\}|,
\]
which is independent of $x\in X$.
The \emph{second eigenmatrix} $Q=(Q_{ij})_{i,j=0}^s$ and the \emph{Krein parameters} $q_{ab}^h$ are specified by
\[
  E_j=\frac1N\sum_{i=0}^sQ_{ij}A_i,
  \qquad
  E_a\circ E_b=\frac1N\sum_{h=0}^sq_{ab}^hE_h
  \qquad(0\leq j,a,b\leq s).
\]
In particular, $Q_{0j}=m_j$; see \cite[Section~2]{MartinTanaka2009} for these parameters.

\begin{proposition}
\label[proposition]{prop:primitive-pair-criterion}
Let $j,k\in\{0,\ldots,s\}$ satisfy $m_j=m_k=3$.
The map $\Phi_{E_j,E_k}$ is injective if and only if
\[
  (Q_{ij},Q_{ik})\notin\{(3,3),(-3,-3)\}
  \qquad(1\leq i\leq s).
\]
If $\Phi_{E_j,E_k}$ is injective, $\Phi_{E_j,E_k}(X)$ is a Grassmannian $4$-design if and only if
\[
  j\neq k
  \qquad\text{and}\qquad
  \sum_{h=0}^sm_h(q_{jk}^h)^2=9.
\]
Moreover,
\begin{equation}
  \sum_{h=0}^sm_h(q_{jk}^h)^2
  =\frac1N\sum_{i=0}^sk_iQ_{ij}^2Q_{ik}^2.
  \label{eq:primitive-pair-sum}
\end{equation}
\end{proposition}

\begin{proof}
For $(x,y)\in R_i$, the spherical embeddings $\varphi_{E_j}$ and $\varphi_{E_k}$ give inner products $Q_{ij}/3$ and $Q_{ik}/3$, respectively.
Thus \cref{lem:faithfulness} gives the injectivity condition.
To check \eqref{eq:second-product-moment}, expand $E_j\circ E_k$ in the orthogonal bases $\{A_i\}_{i=0}^s$ and $\{E_h\}_{h=0}^s$ of $\mathcal A$.
Since $\lVert A_i\rVert_F^2=Nk_i$ and $\lVert E_h\rVert_F^2=m_h$, this gives
\[
  N^2\lVert E_j\circ E_k\rVert_F^2
  =\frac1N\sum_{i=0}^sk_iQ_{ij}^2Q_{ik}^2
  =\sum_{h=0}^sm_h(q_{jk}^h)^2.
\]
The Gram matrices are $NE_j/3$ and $NE_k/3$, so their entrywise product has squared Frobenius norm $\frac{N^2}{81}\sum_{h=0}^sm_h(q_{jk}^h)^2$.
Equation \eqref{eq:second-product-moment} is therefore equivalent to $\sum_{h=0}^sm_h(q_{jk}^h)^2=9$.
Equation \eqref{eq:scheme-first-moment} shows that \eqref{eq:first-moment} holds exactly when $j\neq k$, since $E_jE_k=0$ for $j\neq k$ and $\tr(E_j^2)=3$.
Conditions \eqref{eq:first-moment} and \eqref{eq:second-product-moment} are necessary and sufficient by \cref{lem:G24-moment-criterion}.
\end{proof}

Even when $\Phi_{E_j,E_k}$ is injective, $\varphi_{E_j}$ need not be faithful.
To apply the classification of faithful spherical embeddings in $\mathbb R^3$ \cite[Theorem~3 in the preprint]{BannaiZhao2020}, we pass to the quotient by the fibers of $\varphi_{E_j}$.

\Needspace{8\baselineskip}
\begin{lemma}
\label[lemma]{lem:coordinate-quotient}
Let $j\in\{0,\ldots,s\}$ satisfy $m_j=3$.
The fibers of $\varphi_{E_j}$ are the classes of the equivalence relation
\[
  \{(x,y)\in X^2:\varphi_{E_j}(x)=\varphi_{E_j}(y)\}
  =
  \bigcup_{\substack{0\leq i\leq s\\Q_{ij}=3}}R_i
\]
and all have size
\[
  f:=\sum_{\substack{0\leq i\leq s\\Q_{ij}=3}}k_i.
\]
Let $Y_j$ be the set of these classes.
The quotient symmetric association scheme on $Y_j$ has a primitive idempotent of rank $3$ with a faithful spherical embedding having image $\varphi_{E_j}(X)$.
In particular, $|Y_j|=N/f$ divides $N$.
\end{lemma}

\begin{proof}
For $(x,y)\in R_i$, we have $\langle\varphi_{E_j}(x),\varphi_{E_j}(y)\rangle=Q_{ij}/3$.
As the vectors have unit length, $Q_{ij}=3$ precisely when $\varphi_{E_j}(x)=\varphi_{E_j}(y)$.
Summing the valencies therefore gives the fiber size $f$.
Since the equivalence relation is a union of relations of the scheme, the distinct relations
\[
  \{(C,D)\in Y_j^2:(C\times D)\cap R_i\neq\varnothing\}
  \qquad(0\leq i\leq s)
\]
form a symmetric quotient association scheme \cite[Section~2, following Lemma~2.5]{MartinTanaka2009}.

Order $X$ class by class. The adjacency matrix of the equivalence relation is
\[
  B:=\sum_{\substack{0\leq i\leq s\\Q_{ij}=3}}A_i
  =I_{|Y_j|}\otimes J_f.
\]
The rows and columns of $E_j$ are constant within each class, so $E_jB=BE_j=fE_j$ and $E_j\in(B/f)\mathcal A(B/f)$.
The Bose--Mesner algebra of the quotient is isomorphic to $(B/f)\mathcal A(B/f)$ by the map $M\mapsto f^{-1}M\otimes J_f$ \cite[Section~2, footnote~12]{MartinTanaka2009}.
Hence there is an idempotent $\overline E_j$ in the quotient algebra such that
\begin{equation}
  E_j=\frac1f\overline E_j\otimes J_f.
  \label{eq:quotient-idempotent}
\end{equation}
Any decomposition of $\overline E_j$ into nonzero orthogonal idempotents in the quotient Bose--Mesner algebra would lift to such a decomposition of $E_j$ in $\mathcal A$, contradicting its primitivity.
Thus $\overline E_j$ is primitive.
Since $J_f$ has rank $1$, \eqref{eq:quotient-idempotent} gives $\rank(\overline E_j)=3$.
For $C,D\in Y_j$, $x\in C$, and $y\in D$, we have
\[
  \frac{|Y_j|}{3}(\overline E_j)_{CD}
  =\frac N3(E_j)_{xy}
  =\langle\varphi_{E_j}(x),\varphi_{E_j}(y)\rangle.
\]
Thus $C\mapsto\varphi_{E_j}(x)$ ($x\in C$) is a spherical embedding of $\overline E_j$ with image $\varphi_{E_j}(X)$. Since the classes are the fibers of $\varphi_{E_j}$, this map is well defined and faithful.
\end{proof}

\begin{lemma}
\label[lemma]{lem:small-order-reduction}
Let $j,k\in\{0,\ldots,s\}$ satisfy $m_j=m_k=3$.
Suppose that $1<N<21$ and $\Phi_{E_j,E_k}$ is injective.
Then $N\in\{4,8,12,16,18,20\}$.
If $\Phi_{E_j,E_k}(X)$ is a Grassmannian $4$-design, then $N\in\{12,16,18,20\}$.
\end{lemma}

\begin{proof}
By \cref{lem:coordinate-quotient} and \cite[Theorem~3 in the preprint]{BannaiZhao2020}, $|Y_j|$ and $|Y_k|$ belong to $\{4,6,8,12,20,30\}$ and divide $N$.
Since $N<21$, it follows that $N\in\{4,6,8,12,16,18,20\}$.
If $N=6$ and $j\neq k$, the orthogonal idempotents $E_0,E_j,E_k$ have total rank $1+3+3>6$, a contradiction.
If $N=6$ and $j=k$, then $|Y_j|=6$, and the image of $\varphi_{E_j}$ is the vertex set of a regular octahedron by \cite[Theorem~3 in the preprint]{BannaiZhao2020}.
Its opposite vertices $v,-v$ give the same plane $[v,v]=[-v,-v]$ under $\Phi_{E_j,E_j}$, contrary to injectivity.
Finally, \cref{lem:G24-moment-criterion} gives $N\geq10$ for a Grassmannian $4$-design, excluding $N=4,8$.
\end{proof}

\subsection{Small-order classification}

We use the complete lists of association schemes in \cite{Nomiyama1995,Hirasaka1997,HanakiMiyamoto1998a,HanakiMiyamoto1998b,HanakiMiyamoto2003}, as implemented in version $3.1.0$ of \texttt{AssociationSchemes} \cite{AssociationSchemes2024}.
For each order $N$, No.~$a$ is the original database index in \texttt{HomogeneousCoherentConfiguration}$(N,a)$; we retain only the symmetric schemes.
The ordering of the idempotents follows the rows of \texttt{MatrixOfEigenvalues}, with $E_j$ corresponding to row $j+1$ \cite[Sections~3.3 and~3.6 of the manual]{AssociationSchemes2024}.

\Needspace{14\baselineskip}
\begin{theorem}
\label[theorem]{thm:small-order-classification}
For symmetric association schemes of order $1<N<21$, the pairs of primitive idempotents $E_j,E_k$ of rank $3$ for which $\Phi_{E_j,E_k}$ is injective are exactly those in the following table, up to exchanging $j$ and $k$.
\begin{center}
\begin{tabular}{@{}cclcc@{}}
\toprule
$N$ & Scheme & $(j,k)$ & $(|Y_j|,|Y_k|)$ & $\sum_{h=0}^sm_h(q_{jk}^h)^2$ \\
\midrule
$4$ & No.~$1$ & $(1,1)$ & $(4,4)$ & $21$ \\
$8$ & No.~$5$ & $(2,3)$ & $(8,4)$ & $21$ \\
$12$ & No.~$20$ & $(2,4),(3,4)$ & $(12,6)$ & $27/2$ \\
$16$ & No.~$12$ & $(1,2)$ & $(4,4)$ & $9$ \\
$16$ & No.~$32$ & $(1,2),(1,3),(2,3)$ & $(4,4)$ & $9$ \\
$16$ & No.~$33$ & $(1,2),(1,3),(2,3)$ & $(4,4)$ & $9$ \\
$16$ & No.~$62$ & $1\leq j<k\leq5$ & $(4,4)$ & $9$ \\
\bottomrule
\end{tabular}
\end{center}
The image of $\Phi_{E_j,E_k}$ is a Grassmannian $4$-design precisely in the four rows with $N=16$.
\end{theorem}

\begin{proof}
By \cref{lem:small-order-reduction}, only the orders $4,8,12,16,18,20$ require inspection.
The numbers of symmetric association schemes and of those having a primitive idempotent of rank $3$ are as follows.
\begin{center}
\begin{tabular}{@{}crrrrrr@{}}
\toprule
$N$ & $4$ & $8$ & $12$ & $16$ & $18$ & $20$ \\
\midrule
Symmetric schemes & $3$ & $10$ & $21$ & $63$ & $32$ & $41$ \\
With a rank-three primitive idempotent & $1$ & $2$ & $7$ & $16$ & $5$ & $5$ \\
\bottomrule
\end{tabular}
\end{center}
For each symmetric scheme in the database, we inspect the pairs
\[
  \{(j,k):0\leq j\leq k\leq s,\ m_j=m_k=3\}.
\]
By \cref{prop:primitive-pair-criterion}, a pair is excluded precisely when $(Q_{ij},Q_{ik})\in\{(3,3),(-3,-3)\}$ for some $i>0$.
For every pair that remains, we calculate $|Y_j|$ and $|Y_k|$ by \cref{lem:coordinate-quotient} and $\sum_{h=0}^sm_h(q_{jk}^h)^2$ by \eqref{eq:primitive-pair-sum}.
These calculations use exact arithmetic in the eigenvalues of the adjacency matrices and give the table in the statement. No pair gives an injective map for $N=18$ or $20$.
By \cref{prop:primitive-pair-criterion}, the image is a Grassmannian $4$-design exactly when $j\neq k$ and $\sum_{h=0}^sm_h(q_{jk}^h)^2=9$, which occurs in the four rows with $N=16$.
\end{proof}

Neither idempotent used in \cref{thm:even-cycle} is primitive, so \cref{thm:small-order-classification} does not apply to that construction.

%%%%%%%%%%%%%%%%%%%%%%%%%%%%%%%%%%%%%%%
\subsection*{Acknowledgements}
The second author is supported by JSPS KAKENHI Grant Number 26KJ0570.
The third author is supported by JSPS KAKENHI Grant Number 25K07095.
The second author used ChatGPT (OpenAI, GPT-5.6 Pro) to check computations.
In particular, the AI identified an error in Lemma~\ref{lem:pair-conic}, which was corrected by the second author. The authors verified the revised proof. The authors take full responsibility for the content of this paper.

\begingroup
\setstretch{1}
\microtypesetup{expansion=false}

\endgroup


\begin{thebibliography}{99}

\bibitem{Bachoc2005}
C. Bachoc, Designs, groups and lattices,
\emph{J. Th\'eor. Nombres Bordeaux}
\textbf{17} (2005), no.~1, 25--44,
\href{https://doi.org/10.5802/jtnb.474}{doi:10.5802/jtnb.474}.

\bibitem{BachocBannaiCoulangeon2004}
C. Bachoc, E. Bannai, and R. Coulangeon,
Codes and designs in Grassmannian spaces,
\emph{Discrete Math.}
\textbf{277} (2004), 15--28,
\href{https://doi.org/10.1016/S0012-365X(03)00151-1}{doi:10.1016/S0012-365X(03)00151-1}.

\bibitem{BachocCoulangeonNebe2002}
C. Bachoc, R. Coulangeon, and G. Nebe,
Designs in Grassmannian spaces and lattices,
\emph{J. Algebraic Combin.}
\textbf{16} (2002), 5--19,
\href{https://doi.org/10.1023/A:1020826329555}{doi:10.1023/A:1020826329555}.

\bibitem{BachocEhler2013}
C. Bachoc and M. Ehler,
Tight $p$-fusion frames,
\emph{Appl. Comput. Harmon. Anal.}
\textbf{35} (2013), 1--15,
\href{https://doi.org/10.1016/j.acha.2012.07.001}{doi:10.1016/j.acha.2012.07.001}.

\bibitem{BannaiIto1984}
E. Bannai and T. Ito,
\emph{Algebraic Combinatorics I: Association Schemes},
Benjamin/Cummings, Menlo Park, CA, 1984.

\bibitem{BannaiZhao2020}
E. Bannai and D. Zhao,
Spherical embeddings of symmetric association schemes in $3$-dimensional Euclidean space,
\emph{Graphs Combin.}
\textbf{36} (2020), 245--250,
\href{https://doi.org/10.1007/s00373-019-02075-6}{doi:10.1007/s00373-019-02075-6};
preprint: \href{https://arxiv.org/abs/1710.00212v1}{arXiv:1710.00212v1}.

\bibitem{AssociationSchemes2024}
J. Bamberg, A. Hanaki, and J. Lansdown,
\texttt{AssociationSchemes}: a GAP package for working with association schemes and homogeneous coherent configurations,
version 3.1.0, 2024,
\href{https://doi.org/10.5281/zenodo.2634954}{doi:10.5281/zenodo.2634954}.

\bibitem{CalderbankCasazzaHeineckeKutyniokPezeshki2011}
R. Calderbank, P. G. Casazza, A. Heinecke, G. Kutyniok, and A. Pezeshki,
Sparse fusion frames: existence and construction,
\emph{Adv. Comput. Math.}
\textbf{35} (2011), no.~1, 1--31,
\href{https://doi.org/10.1007/s10444-010-9162-3}{doi:10.1007/s10444-010-9162-3}.

\bibitem{CalderbankHardinRainsShorSloane1999}
A. R. Calderbank, R. H. Hardin, E. M. Rains, P. W. Shor, and N. J. A. Sloane,
A group-theoretic framework for the construction of packings in Grassmannian spaces,
\emph{J. Algebraic Combin.}
\textbf{9} (1999), no.~2, 129--140,
\href{https://doi.org/10.1023/A:1018673825179}{doi:10.1023/A:1018673825179}.

\bibitem{ConwayHardinSloane1996}
J. H. Conway, R. H. Hardin, and N. J. A. Sloane,
Packing lines, planes, etc.: packings in Grassmannian spaces,
\emph{Experiment. Math.}
\textbf{5} (1996), 139--159,
\href{https://doi.org/10.1080/10586458.1996.10504585}{doi:10.1080/10586458.1996.10504585}.

\bibitem{deLaHarpePache2005}
P. de la Harpe and C. Pache,
Cubature formulas, geometrical designs, reproducing kernels, and Markov operators,
in \emph{Infinite Groups: Geometric, Combinatorial and Dynamical Aspects},
Progr. Math., vol.~248, Birkh\"auser, Basel, 2005, pp.~219--267,
\href{https://doi.org/10.1007/3-7643-7447-0_6}{doi:10.1007/3-7643-7447-0\_6}.

\bibitem{DelsarteGoethalsSeidel1977}
P. Delsarte, J. M. Goethals, and J. J. Seidel,
Spherical codes and designs,
\emph{Geom. Dedicata}
\textbf{6} (1977), no.~3, 363--388,
\href{https://doi.org/10.1007/BF03187604}{doi:10.1007/BF03187604}.

\bibitem{EtayoLopezGomez2025}
U. Etayo and P. R. L\'opez-G\'omez,
Minimal Riesz and logarithmic energies on the Grassmannian $\operatorname{Gr}_{2,4}$,
preprint,
\href{https://arxiv.org/abs/2501.00424}{arXiv:2501.00424}.

\bibitem{HanakiMiyamoto1998a}
A. Hanaki and I. Miyamoto,
Classification of association schemes with 16 and 17 vertices,
\emph{Kyushu J. Math.}
\textbf{52} (1998), 383--395,
\href{https://doi.org/10.2206/kyushujm.52.383}{doi:10.2206/kyushujm.52.383}.

\bibitem{HanakiMiyamoto1998b}
A. Hanaki and I. Miyamoto,
Classification of association schemes with 18 and 19 vertices,
\emph{Korean J. Comput. Appl. Math.}
\textbf{5} (1998), 543--551,
\href{https://doi.org/10.1007/BF03008881}{doi:10.1007/BF03008881}.

\bibitem{HanakiMiyamoto2003}
A. Hanaki and I. Miyamoto,
Classification of association schemes of small order,
\emph{Discrete Math.}
\textbf{264} (2003), 75--80,
\href{https://doi.org/10.1016/S0012-365X(02)00551-4}{doi:10.1016/S0012-365X(02)00551-4}.

\bibitem{HanakiMiyamotoCatalogue}
A. Hanaki and I. Miyamoto,
Classification of association schemes, online catalogue,
\url{http://math.shinshu-u.ac.jp/~hanaki/as/};
data file
\href{http://math.shinshu-u.ac.jp/~hanaki/as/data/as10}{\texttt{as10}}.

\bibitem{Hirasaka1997}
M. Hirasaka,
The classification of association schemes with $11$ or $12$ vertices,
\emph{Kyushu J. Math.}
\textbf{51} (1997), 413--428,
\href{https://doi.org/10.2206/kyushujm.51.413}{doi:10.2206/kyushujm.51.413}.

\bibitem{LansdownMartinMunemasaSudaTanaka2026}
J. Lansdown, W. J. Martin, A. Munemasa, S. Suda, and H. Tanaka,
On the spherical design properties of a $P$- and $Q$-polynomial association scheme,
preprint, 2026,
\href{https://arxiv.org/abs/2608.27082}{arXiv:2608.27082}.

\bibitem{MartinTanaka2009}
W. J. Martin and H. Tanaka,
Commutative association schemes,
\emph{European J. Combin.}
\textbf{30} (2009), no.~6, 1497--1525,
\href{https://doi.org/10.1016/j.ejc.2008.11.001}{doi:10.1016/j.ejc.2008.11.001}.

\bibitem{Misawa2026}
R. Misawa,
Constructing spherical designs using tight $t$-fusion frames,
preprint, 2026,
\href{https://arxiv.org/abs/2601.17294v1}{arXiv:2601.17294v1}.

\bibitem{Nomiyama1995}
E. Nomiyama,
Classification of association schemes with at most ten vertices,
\emph{Kyushu J. Math.}
\textbf{49} (1995), 163--195,
\href{https://doi.org/10.2206/kyushujm.49.163}{doi:10.2206/kyushujm.49.163}.

\bibitem{Suda2011}
S. Suda,
On spherical designs obtained from $Q$-polynomial association schemes,
\emph{J. Combin. Des.}
\textbf{19} (2011), no.~3, 167--177,
\href{https://doi.org/10.1002/jcd.20278}{doi:10.1002/jcd.20278}.

\bibitem{Weyl1946}
H. Weyl,
\emph{The Classical Groups: Their Invariants and Representations},
2nd ed., Princeton Mathematical Series, vol.~1,
Princeton University Press, Princeton, NJ, 1946.

\end{thebibliography}
\end{document}